\mag=\magstep1
\documentstyle[draft,amssymb]{amsart}
\newtheorem{thm}{Theorem}[subsection]
\newtheorem{lem}[thm]{Lemma}
\newtheorem{prop}[thm]{Proposition}
\newtheorem{conj}[thm]{Conjecture}

\theoremstyle{definition}
\newtheorem{defn}[thm]{Definition}

\theoremstyle{remark}
\newtheorem{rem}[thm]{Remark}

\renewcommand{\theequation}{\thesubsection.\arabic{equation}}

\newcommand{\thmref}[1]{Theorem~\ref{#1}}

\newcommand{\secref}[1]{\S\ref{#1}}
\newcommand{\lemref}[1]{Lemma~\ref{#1}}
\newcommand{\appref}[1]{Appendix~\ref{#1}}
\newcommand{\eps}{\varepsilon}
\newcommand{\lam}{\lambda} \newcommand{\Lam}{\Lambda}
\newcommand{\sig}{\sigma}
\newcommand{\calH}{{\cal H}}
\newcommand{\bbbS}{{\Bbb S}}

\newcommand{\calT}{{\cal T}}
\newcommand{\Integer}{{\Bbb Z}}
\newcommand{\Real}{{\Bbb R}}
\newcommand{\Complex}{{\Bbb C}}
\newcommand{\End}{\operatorname{End}}
\newcommand{\ord}{\operatorname{ord}}
\newcommand{\Tr}{\operatorname{Tr}}
\newcommand{\tensor}{\otimes}
\renewcommand{\Re}{\operatorname{Re}}
\renewcommand{\Im}{\operatorname{Im}}
\newcommand{\sh}{\operatorname{sinh}}
\newcommand{\ch}{\operatorname{cosh}}
\newcommand{\half}{{\textstyle \frac{1}{2}}}
\newcommand{\faceW}[5]{W\left(\left.\begin{matrix}#1&#2\\#3&#4%
\end{matrix}\right|#5\right)}
\begin{document}

\title[Higher Spin 8 Vertex Model]
{Bethe Ansatz for higher spin\\
eight vertex models}

\author{Takashi Takebe}
\address{Department of Mathematical Sciences\\
        the University of Tokyo, Hongo 7-3-1, Bunkyo-ku\\
        Tokyo, 113 Japan}
\email{takebe@@math.s.u-tokyo.ac.jp}


\dedicatory{Dedicated to the memory of Alexander A.~Belov}


\maketitle
\centerline{Department of Mathematical Sciences,}
\centerline{the University of Tokyo,}
\centerline{Hongo 7-3-1, Bunkyo-ku, Tokyo, 113 Japan}
\overfullrule=0pt

\begin{abstract}
A generalization of the eight vertex model by means of higher spin
representations of the Sklyanin algebra is investigated by the quantum
inverse scattering method and the algebraic Bethe Ansatz. Under the
well-known string hypothesis low-lying excited states are considered
and scattering phase shifts of two physical particles are
calculated. The $S$ matrix of two particle states is shown to be
proportional to the Baxter's elliptic $R$ matrix with a different
elliptic modulus from the original one.
\end{abstract}
%
%
%
%
%
%
\section*{Introduction}
\label{intro}
In this paper we consider a generalization of the eight vertex model
by means of higher spin representations (spin $\ell$) of the Sklyanin
algebra \cite{skl:alg} on a space of theta functions \cite{skl:rep}.
This model has $2\ell+1$ dimensional state space on each vertical edge
and two dimensional state space on each horizontal edge.

Relation of the eight vertex model to the SOS type model was
established by Baxter \cite{bax:eight}. A similar relation also holds
in our case and, using this relation, we can pursue the quantum inverse
scattering method and the algebraic Bethe Ansatz, following
\cite{kul-skl}, \cite{takh-fad:xyz}, \cite{takh-fad:xxx}. In the first
part of this paper, we examine Bethe vectors and give coordinate
expression of them in terms of Boltzmann weights of SOS type model.
We also prove a sum rule of rapidities of quasi-particles, which was
proved for the eight vertex model by Baxter \cite{bax:part} who made
use of a functional equation, an alternative to Bethe Ansatz. This
rule is related to a parity of Bethe vectors.

A higher spin version of SOS type model was constructed by Date,
Jimbo, Miwa and Okado \cite{djmo}, \cite{jmo}, \cite{djkmo} by fusion
procedure \cite{kul-resh-skl}, \cite{cher}, \cite{zhou-hou},
\cite{hou-zhou}. Recently Hasegawa showed \cite{has:prep} that
a representation of the Sklyanin algebra obtained by fusion procedure
repeated $2\ell-1$ times is equivalent to the spin $\ell$
representation on a space of theta functions \cite{skl:rep}.
Hence, in principle our model is equivalent to the higher spin SOS
model by Date and others.
The use of representations by Sklyanin \cite{skl:rep} makes it possible
to compute eigenvectors of transfer matrices explicitly and to apply
the quantum inverse scattering method and the algebraic Bethe Ansatz
directly (cf.~\cite{bazh-resh}).

As is shown by Baxter \cite{bax:xyz}, the transfer matrix of the eight
vertex model contains the Hamiltonian of an anisotropic Heisenberg
magnetic chain (the $XYZ$ model) \cite{hei}.
Our model is also related to a quantum spin chain model
with $2\ell+1$ dimensional local quantum spaces.
But in general the transfer matrix of our lattice model does not give
local Hamiltonians directly, because the dimension of the auxiliary
space is fixed to two. In order to write down the Hamiltonian
of this spin chain, we must use the fused transfer matrix
corresponding to a $2\ell+1$ dimensional auxiliary space.
More generally we can construct a model with arbitrary spins on quantum
and auxiliary spaces by fusion procedure.
We will study such models in the forthcoming paper.

Note that Bethe vectors constructed above give eigenvectors of these models
simultaneously under assumption of non-degeneracy,
since transfer matrices with different auxiliary spins are mutually
commuting \cite{kul-resh-skl}.

Though momenta and Hamiltonians of spin chains are calculated from
transfer matrices of fused models, $S$ matrices (phase shifts) of spin
waves do not depend on the auxiliary space. In the second part of this
paper, we calculate a two particle $S$ matrix of spin waves from Bethe
vectors obtained above, following the recipe by Korepin \cite{kor},
Destri and Lowenstein \cite{des-low}.
The result confirms Smirnov's conjecture \cite{fijkmy} which states that
this $S$ matrix should be given by an elliptic $R$ matrix, the
elliptic modulus of which is different from that of the original $R$
matrix in the definition of the model.

Corresponding results were established for the totally isotropic
models (the $XXX$ model) and its higher spin generalization, by
Faddeev, Takhtajan, Babujian, Avdeev and D\"orfel \cite{takh-fad:xxx},
\cite{takh}, \cite{bab}, \cite{av-dor} and for the $XXZ$ model
(the six vertex model) and its higher spin generalization, by Sogo,
Kirillov, Reshetikhin \cite{sog}, \cite{kir-resh}. Free energy of
the eight vertex model was obtained by Baxter \cite{bax:part} and
low-lying excited states were studied by Johnson, Krinsky and McCoy
\cite{jo-kr-mc}, but our calculation of the $S$ matrix seems to be new
even for the eight vertex model ($\ell = 1/2$), though a partial
result on the $S$ matrix for this case was calculated by Freund
and Zabrodin \cite{fre-zab}. The algebraic Bethe
Ansatz was shown to be applicable to the higher spin eight vertex
models in \cite{take1} and their free energy was calculated in
\cite{take2}.

This paper is organized as follows:
In Chapter I we begin with review of definition of the model and
generalized algebraic Bethe Ansatz, following \cite{take1}. Then,
giving a whole set of intertwining vectors explicitly, we write down
coordinate expression of Bethe vectors in terms of them. A sum rule of
rapidities of quasi-particles are presented which helps solving the
Bethe equations. The proof is given in \appref{sumrule:pf}.
In Chapter II we study thermodynamic limit of several Bethe vectors.
Free energy is calculated in \secref{gr.state} (this result was
announced in \cite{take2}) and low-lying excited states are examined
in \secref{s-matrix} under assumptions of string configuations. In
particular we compute a two particle $S$ matrix in \secref{s-matrix}.
We summarize prerequisites on the Sklyanin algebra and its
representations in \appref{skl-alg}.

Chapter I and Appendices are of algebraic nature, while in Chapter II
we do not try to give mathematically rigorous arguments to every
detail, since the goal of this chapter is to compute quantities of
physical importance. In order to make this computation rigorous hard
analysis is indispensable, which is beyond the scope of this paper.

%
%
%
%
\section{Description of the model and algebraic Bethe Ansatz}

The model was defined and its eigenvectors were constructed by the
generalized algebraic Bethe Ansatz in \cite{take1}. In
\secref{model-def} and \secref{aBA} we briefly review this work since
we change notations and normalizations. In \secref{intertwine} we find
an explicit form of a whole set of intertwining vectors which enables
us to write down coordinate expression of Bethe vectors
(\secref{aBA}). Only the ``highest'' ones of intertwining vectors
(``local pseudo vacua'' in the context of the algebraic Bethe Ansatz
\cite{takh-fad:xyz}) were used in \cite{take1} to construct Bethe
vectors. In \secref{sumrule} we show that sum of rapidities of
quasi-particles should satisfy an integrality condition. This comes
from quasi-periodicity of theta functions and therefore is absent in
the case of models associated to trigonometric and rational $R$ matrices.

\subsection{Definition of the model}
\setcounter{equation}{0}
\label{model-def}

The model is parametrized by a half integer $\ell$ and two complex
parameters: an elliptic modulus $\tau$ and an anisotropy parameter
$\eta$. In this paper we assume that the elliptic modulus is a pure
imaginary number while the anisotropy parameter is a rational number:
\begin{equation}
    \tau = \frac{i}{t}, \qquad \eta = \frac{r'}{r},
\end{equation}
where $t>0$ and  $r$, $r'$ are integers mutually coprime. Moreover we
impose a condition that $r$ is even, $r'$ is odd, and
$2(2\ell+1) \eta < 1$.

Now we define a lattice model of vertex type as in \cite{take1}.
We consider a square lattice with $N$ columns and $N'$ rows on a torus,
i.e., periodic boundary condition imposed.
States on the $n$-th vertical edge belong to the spin $\ell$ representation
space $V_n \simeq \Theta^{4\ell+}_{00}$
of the Sklyanin algebra (see \appref{skl-alg}) while states on each
horizontal edge are two dimensional vectors. A row-to-row transfer
matrix, $T(\lam)$, of the model is defined as the trace of a monodromy
matrix, $\calT(\lam)$, in the context of the quantum inverse
scattering method \cite{takh-fad:xyz}:
\begin{gather}
    \calT(\lam)
    = \begin{pmatrix} A_N(\lam) & B_N(\lam) \\
                      C_N(\lam) & D_N(\lam) \end{pmatrix}
    := L_N(\lam) \dots L_2(\lam) L_1(\lam),
\label{def:monodromy}
\\
    T(\lam)     = \Tr_{\Complex^2} (\calT(\lam)) = A_N(\lam) + D_N(\lam),
\label{def:transfer}
\end{gather}
where the $L$ operators, $L_n(\lam)$, are defined by (cf.~\eqref{def:L})
\begin{equation}
\begin{gathered}
    L_n(\lam)
    = \begin{pmatrix} \alpha_n(\lam) & \beta_n(\lam)   \\
                      \gamma_n(\lam) & \delta_n(\lam) \end{pmatrix}
    := \sum_{a=0}^3 W_a^L(\lam) \rho^\ell_n(S^a) \tensor \sig^a,\\
    \rho^\ell_n(S^a) = 1 \tensor \dots \tensor 1
                         \tensor \rho^\ell(S^a) \tensor
                       1 \tensor \dots \tensor 1,
\end{gathered}
\label{def:L-n}
\end{equation}
elements of which act on a Hilbert space
$\calH = \bigotimes_{n=1}^N V_n$, but non-trivially
only on the $n$-th component.
Assignment of Boltzmann weights to vertices are determined by this
$L$ operator. The monodromy matrix is a $2\times2$ matrix with
elements in $\End_{\Complex}(\calH)$, and the transfer matrix is an
operator in $\End_{\Complex}(\calH)$.

The partition function, $Z(\lam)$, and the free energy per site,
$f(\lam)$, are
\begin{align*}
    Z(\lam) &= \Tr_{\calH} (T(\lam)^{N'}),
\\
    - \beta f(\lam) &= \frac{1}{NN'} \log Z(\lam).
\end{align*}
In the thermodynamic limit, $N,N' \to\infty$, only the greatest
eigenvalue, $\Lam_{\text{max}}$, of $T(\lam)$ contribute to the free
energy
\begin{equation}
    - \beta f(\lam) =
    \lim_{N \to \infty} \frac{1}{N} |\Lam_{\text{max}}|,
\label{free-energy}
\end{equation}
which was computed in \cite{take2}. We will recall this result in
\secref{gr.state} with details omitted in \cite{take2}.

\begin{rem}
The above defined model is a {\em homogeneous} lattice in the sense
that it is invariant with respect to vertical and horizontal translation.
We can also define an {\em inhomogeneous} lattice by assigning different
spectral parameter, $\lam_i$, and different spin, $\ell_i$, to each
vertical edge. We have only to replace $L_n(\lam)$ in \eqref{def:monodromy}
and \eqref{def:transfer} by
\begin{gather*}
    L_n(\lam-\lam_n)
    := \sum_{a=0}^3
       W_a^L(\lam-\lam_n) \rho^{\ell_n}_n(S^a) \tensor \sig^a,\\
    \rho^{\ell_n}_n (S^a) = 1 \tensor \dots \tensor 1
                              \tensor \rho^{\ell_n}(S^a) \tensor
                            1 \tensor \dots \tensor 1.
\end{gather*}
All arguments in \secref{intertwine}, \secref{aBA} remain true
with suitable changes, as is shown in \cite{take1}. Such models are
important for the study of certain integrable systems \cite{kuz}.
\end{rem}

\subsection{Intertwining vectors, gauge transformation}
\setcounter{equation}{0}
\label{intertwine}

Intertwining vectors were first introduced by Baxter \cite{bax:eight},
and given an interpretation as a gauge transformation in the context
of the quantum inverse scattering method by Takhtajan and Faddeev
\cite{takh-fad:xyz}. Generalization to higher spin case by means of
fusion procedure was studied by Date, Jimbo, Miwa, Kuniba and Okado
\cite{djmo}, \cite{djkmo}, \cite{jmo}.
Here we define intertwining vectors directly in the space of theta
fucntions. They should be identified with those defined in \cite{djmo},
\cite{djkmo}, \cite{jmo} through Hasegawa's isomorphism
\cite{has:prep}.

\begin{defn}
Let $k, k'$ be integers satisfying
$k - k' \in \{ -2\ell, -2\ell + 2, \dots, 2\ell - 2, 2\ell\}$,
and $\lam$, $s=(s_+, s_-)$ be complex parameters.
We call the following vectors
$\phi_{k,k'}(\lam;s) = \phi^{(\ell)}_{k,k'} (\lam;s)(z)
\in \Theta^{4\ell+}_{00}$ {\em intertwining vectors} of spin $\ell$:
\begin{equation}
\begin{aligned}
    \phi^{(\ell)}_{k,k'} (\lam;s)(z) = a_{k,k'}
    \prod_{j=1}^{\ell + \frac{k-k'}{2}}
    & \theta\left(
            z + \frac{s_+ - \lam}{2} + \frac{\tau}{4} + (k'-\ell+2j-1)\eta
            \right) \times \\
    & \times
      \theta\left(
            z - \frac{s_+ - \lam}{2} - \frac{\tau}{4} - (k'-\ell+2j-1)\eta
            \right) \\
    \prod_{j=1}^{\ell - \frac{k-k'}{2}}
    & \theta\left(
            z + \frac{s_- + \lam}{2} + \frac{\tau}{4} + (k -\ell+2j-1)\eta
            \right) \times \\
    & \times
      \theta\left(
            z - \frac{s_- + \lam}{2} - \frac{\tau}{4} - (k -\ell+2j-1)\eta
            \right).
\end{aligned}
\label{def:int-vec}
\end{equation}
Here $\theta(z) = \theta_{00}(z;\tau)$ (see \eqref{theta:def}), and
$a_{k,k'} = e^{2\pi i \ell (k+k') \eta}
            (i t^{-1/2} e^{\pi i (s_+-s_-)})^{\frac{k-k'}{2}}$.
\end{defn}

Following \cite{takh-fad:xyz}, we introduce a matrix of {\em gauge
transformation} $M_k$:
\begin{equation}
\begin{split}
    M_k (\lam;s) =&
    \begin{pmatrix}
    \theta_{11}(-it(s_+ - \lam + 2k\eta);2it) &
    \theta_{11}(-it(s_- + \lam + 2k\eta);2it) \\
    \theta_{01}(-it(s_+ - \lam + 2k\eta);2it) &
    \theta_{01}(-it(s_- + \lam + 2k\eta);2it) \\
    \end{pmatrix} \times
    \\
    \times &
    \begin{pmatrix}
    e^{-\frac{\pi t}{2} ( s_+ - \lam+2k\eta-\frac{i}{2t} )^2}&
    0 \\
    0 &
    e^{-\frac{\pi t}{2} ( s_- + \lam+2k\eta-\frac{i}{2t} )^2}
    \end{pmatrix} \times
    \\
    \times &
    \begin{pmatrix}
    1 & 0 \\ 0 & \theta_{11}(w_k;\tau)^{-1}
    \end{pmatrix}
\end{split}
\label{def:M-k}
\end{equation}
where
\begin{equation}
    w_k = \frac{s_+ + s_-}{2} + 2 k \eta - \frac{\tau}{2}.
\label{def:w-k}
\end{equation}

Let us define a twisted $L$ operator by
\begin{equation}
\begin{split}
    L_{k,k'}(\lam;s) &=
    \begin{pmatrix}
    \alpha_{k,k'}(\lam;s) & \beta_{k,k'} (\lam;s) \\
    \gamma_{k,k'}(\lam;s) & \delta_{k,k'}(\lam;s)
    \end{pmatrix} \\
    &:= M_k^{-1}(\lam;s) L(\lam) M_{k'}(\lam;s).
\end{split}
\end{equation}

\begin{prop}
\label{prop:vertex-face}
Each component of $L_{k,k'}$ acts on the intertwining vector as follows:
\begin{equation}
\begin{aligned}
    \alpha_{k,k'}(\lam;s) \phi_{k,k'} &=
    \faceW{k}{k'}{k-1}{k'-1}{\lam}
    \phi_{k-1,k'-1} ,
\\
    \beta_{k,k'}(\lam;s) \phi_{k,k'} &=
    \faceW{k}{k'}{k-1}{k'+1}{\lam}
    \phi_{k-1,k'+1} ,
\\
    \gamma_{k,k'}(\lam;s) \phi_{k,k'} &=
    \faceW{k}{k'}{k+1}{k'-1}{\lam}
    \phi_{k+1,k'-1} ,
\\
    \delta_{k,k'}(\lam;s) \phi_{k,k'} &=
    \faceW{k}{k'}{k+1}{k'+1}{\lam}
    \phi_{k+1,k'+1},
\end{aligned}
\label{vertex-face}
\end{equation}
where $\phi_{k,k'}=\phi^{(\ell)}_{k,k'}(0;s)$ and $W$ is the Boltzmann
weight of SOS type \cite{bax:eight}, \cite{djmo}:
\begin{equation}
\begin{aligned}
    \faceW{k}{k'}{k-1}{k'-1}{\lam} &=
    2 \theta_{11}(\lam + (k-k')\eta)
    \frac{\theta_{11}(w_{(k+k' + 2\ell)/2} )}
         {\theta_{11}(w_k)},
\\
    \faceW{k}{k'}{k-1}{k'+1}{\lam} &=
    2 \theta_{11}((k'-k-2\ell)\eta)
    \frac{\theta_{11}(w_{(k+k')/2} + \lam )}
         {\theta_{11}(w_k) \theta_{11}(w_{k'})},
\\
    \faceW{k}{k'}{k+1}{k'-1}{\lam} &=
    2 \theta_{11}((k-k'-2\ell)\eta)
    \theta_{11}(w_{(k+k')/2} - \lam),
\\
    \faceW{k}{k'}{k+1}{k'+1}{\lam} &=
    2 \theta_{11}(\lam - (k-k')\eta)
    \frac{\theta_{11}(w_{(k+k' - 2\ell)/2})}
         {\theta_{11}(w_k')}.
\end{aligned}
\label{face-weight}
\end{equation}
\end{prop}

This proposition is proved in the same way as (3.7) of \cite{take1}.

If we denote components of $L_{k,k'}$ by
\begin{equation}
    L_{k,k'}(\lam;s) =
    \begin{pmatrix}
    L_{k,k'}(-1,-1;\lam;s) & L_{k,k'}(-1,+1;\lam;s) \\
    L_{k,k'}(+1,-1;\lam;s) & L_{k,k'}(+1,+1;\lam;s)
    \end{pmatrix},
\end{equation}
\eqref{vertex-face} takes the form:
\begin{equation}
    L_{k,k'}(\eps,\eps';\lam;s) \phi_{k,k'} =
    \faceW{k}{k'}{k+\eps}{k'+\eps'}{\lam} \phi_{k+\eps,k'+\eps'}.
\label{vertex-face'}
\end{equation}

In \cite{take1}, vectors
$\omega^n_m = \phi_{n+2\ell m, n + 2\ell(m-1)}(s)$ were called local
vacua. For these vectors the formulae \eqref{vertex-face} reduce to
\begin{align}
    \alpha_{k,k-2\ell} \phi_{k,k-2\ell} (s) &=
    2 \theta_{11}(\lam + 2 \ell \eta) \phi_{k-1,k-2\ell-1} (s),
\label{alpha:vac}
\\
    \delta_{k,k-2\ell} \phi_{k,k-2\ell} (s) &=
    2 \theta_{11}(\lam - 2 \ell \eta) \phi_{k+1,k-2\ell+1} (s),
\label{delta:vac}
\\
    \gamma_{k,k-2\ell} \phi_{k,k-2\ell} (s) &= 0.
\label{gamma:vac}
\end{align}
This property is important for algebraic Bethe Ansatz.

\begin{rem}
Denoting the column vectors of $M_k$ by $\psi_{k,k\pm1}(\lam;s)$,
$M_k=(\psi_{k,k-1}(\lam;s), \psi_{k,k+1}(\lam;s))$,
one can rewrite \eqref{vertex-face} as follows:
\begin{multline*}
    L(\lam-\mu) \phi^{(\ell)}_{k,k'} (\lam;s)\tensor
                \psi_{k',k'+\eps'}(\mu;s)
    = \\ =
    \sum_{\eps} \faceW{k}{k'}{k+\eps}{k'+\eps'}{\lam-\mu}
                \phi^{(\ell)}_{k+\eps,k'+\eps'}(\lam;s) \tensor
                \psi_{k,k+\eps}(\mu;s).
\end{multline*}
Namely $\phi$ and $\psi$ intertwine the vertex weights and the SOS
weights. This is where the name ``intertwining vector'' comes from.
See \cite{bax:eight}, \cite{djmo}, \cite{has:cross}.
Note that $\phi^{(1/2)}_{k, k\pm1}(\lam;s)$ are proportional
to column vectors of $M_k$ under the identification
\eqref{rep:pauli}.
\end{rem}

\subsection{Generalized algebraic Bethe Ansatz}
\setcounter{equation}{0}
\label{aBA}
In this section we recall the construction of eigenvectors of the
transfer matrix by means of the algebraic Bethe Ansatz,
following \cite{take1}, and give several properties of them.
Hereafter we assume that $M:= N\ell$ is an integer.

First introduce a modified monodromy matrix twisted by gauge
transformation:
\begin{equation}
\begin{split}
    \calT_{k,k'}(\lam;s) &=
    \begin{pmatrix}
    A_{k,k'}(\lam;s) & B_{k,k'} (\lam;s) \\
    C_{k,k'}(\lam;s) & D_{k,k'}(\lam;s)
    \end{pmatrix}
\\
    &:= M_k^{-1}(\lam;s) \calT(\lam) M_{k'}(\lam;s),
\end{split}
\label{def:twist-monod}
\end{equation}
and {\em fundamental vectors} in $\calH$:
\begin{equation}
    |a_N, a_{N-1}, \dots, a_1, a_0 \rangle :=
    \phi_{a_N,a_{N-1}} \tensor \dots \tensor
    \phi_{a_2,a_1} \tensor \phi_{a_1,a_0},
\end{equation}
where $\phi_{a,b} = \phi_{a,b}(0;s)$ are intertwining vectors defined
by \eqref{def:int-vec}. We fix a value of the parameter $s=(s_+,s_-)$ and
suppress it unless it is necessary. A {\em pseudo vacuum}, $\Omega^a_N$,
is a fundamental vector characterized by $a_0 = a$,
$a_i - a_{i-1} = 2\ell$ for all $i = 1, \dots, N$:
\begin{equation}
    \Omega^a_N = |a+2N\ell, a+2(N-1)\ell, \dots, a+2\ell, a \rangle.
\end{equation}
This vector satisfies
\begin{align}
    A_{a+2N\ell,a}(\lam) \Omega^a_N &=
    (2 \theta_{11}(\lam + 2\ell\eta))^N \Omega^{a-1}_N,
\label{A:vac}
\\
    D_{a+2N\ell,a}(\lam) \Omega^a_N &=
    (2 \theta_{11}(\lam - 2\ell\eta))^N \Omega^{a+1}_N,
\label{D:vac}
\\
    C_{a+2N\ell,a}(\lam) \Omega^a_N &= 0
\label{C:vac}
\end{align}
by virtue of \eqref{alpha:vac}, \eqref{delta:vac}, \eqref{gamma:vac},
respectively.

As is shown in \cite{take1}, the {\em algebraic Bethe Ansatz} for our
case leads to the following:
\begin{prop}
Let $\nu$ be an integer, $\lam_1, \dots, \lam_M$ complex numbers. Define
a vector $\Psi_\nu(\lam_1, \dots, \lam_M) \in \calH$ by
\begin{equation}
\begin{aligned}
    \Psi_\nu(\lam_1, \dots, \lam_M) &:=
    \sum_{a=0}^{r-1} e^{2\pi i \nu \eta a}
    \Phi_a(\lam_1, \dots, \lam_M),
\\
    \Phi_a(\lam_1, \dots, \lam_M) &:=
    B_{a+1,a-1}(\lam_1) B_{a+2,a-2}(\lam_2) \dots \\
    &\dots B_{a+M,a-M}(\lam_M) \Omega^{a-M}_N.
\end{aligned}
\label{Bethe-vec}
\end{equation}
Then $\Psi_\nu(\lam_1, \dots, \lam_M)$ is an eigenvector of the transfer
matrix $T(\lam)$ with an eigenvalue
\begin{equation}
\begin{aligned}
    t(\lam) &=
    e^{ 2\pi i \nu \eta}
    (2\theta_{11}(\lam + 2\ell\eta;\tau))^N
    \prod_{j=1}^M
    \frac{\theta_{11}(\lam - \lam_j - 2\eta;\tau)}
         {\theta_{11}(\lam - \lam_j;       \tau)} \\
    &+
    e^{-2\pi i \nu \eta}
    (2\theta_{11}(\lam - 2\ell\eta;\tau))^N
    \prod_{j=1}^M
    \frac{\theta_{11}(\lam - \lam_j + 2\eta;\tau)}
         {\theta_{11}(\lam - \lam_j;       \tau)},
\end{aligned}
\label{Bethe-ev}
\end{equation}
provided that $\nu$ and $\{\lam_1, \dots, \lam_M\}$ satisfy the
following {\em Bethe equations}:
\begin{equation}
    \left(\frac{\theta_{11}(\lam_j + 2\ell\eta;\tau)}
               {\theta_{11}(\lam_j - 2\ell\eta;\tau)}\right)^N
    =
    e^{-4\pi i \nu \eta}
    \prod\begin{Sb}k=1 \\ k \neq j \end{Sb}^M
    \frac{\theta_{11}(\lam_j - \lam_k + 2\eta;\tau)}
         {\theta_{11}(\lam_j - \lam_k - 2\eta;\tau)},
\label{Bethe-eq}
\end{equation}
for all $j=1, \dots, M$.
\end{prop}
\begin{pf}
The proof is the same as that in \cite{takh-fad:xyz}. We only recall the
periodicity of the vector $\Phi_a$ with respect to $a$, which is the
reason that we do not have to take an infinite sum in the definition of
$\Psi_\nu$.

Recall that $\eta$ is a rational number $r'/r$. Therefore $2(k+r)\eta
\equiv 2k\eta$ $\pmod{2}$. This fact and quasi periodicity of theta
functions imply $M_{k+r} = M_k$. Hence (see \eqref{def:twist-monod})
\begin{equation}
    \calT_{k+r, k'+r}(\lam) = \calT_{k,k'}(\lam),
\end{equation}
in particular, $B_{k+r, k'+r}(\lam) = B_{k,k'}(\lam)$.
Similarly one can prove $\phi_{k+r,k'+r} = \phi_{k,k'}$. This proves
$\Phi_{a+r}(\lam_1, \dots, \lam_M) = \Phi_a(\lam_1, \dots, \lam_M)$.
\end{pf}

The eigenvalue \eqref{Bethe-ev} is written in a compact form in terms
of a function $Q(\lam)$ defined by
\begin{equation}
    Q(\lam) = e^{-\pi i \nu \lam}
    \prod_{j=1}^M \theta_{11}(\lam - \lam_j).
\label{def:Q}
\end{equation}
The eigenvalue of the transfer matrix for a Bethe vector
$\Psi_\nu(\lam_1, \dots,\lam_M)$ is
\begin{equation}
    t(\lam):= h(\lam+2\ell\eta) \frac{Q(\lam-2\eta)}{Q(\lam)}
            + h(\lam-2\ell\eta) \frac{Q(\lam+2\eta)}{Q(\lam)},
\label{Bethe-ev:Q}
\end{equation}
where $h(z) = (2\theta_{11}(z))^N$.
The Bethe equations \eqref{Bethe-eq} can be interpreted
as the condition of cancellation of poles at $\lam_j$ of the right
hand side of the above equation. This observation is due to Baxter
\cite{bax:part} and a starting point of Reshetikhin's analytic Bethe
Ansatz \cite{resh:analBA}. We essentially use this observation
to derive the sum rule in \secref{sumrule}.

Because of the commutation relation
$
    B_{k,k'+1}(\lam) B_{k+1,k'}(\mu) =
    B_{k,k'+1}(\mu) B_{k+1,k'}(\lam),
$
which is a consequence of relation \eqref{RLL},
$\Psi_\nu(\lam_1, \dots, \lam_M)$ does not depend on the order of parameters
$\lam_1, \dots, \lam_M$. Moreover, we may restrict the parameters to
the fundamental domain
\begin{equation}
    - \frac{1}{2} \leqq \Re \lam_j \leqq \frac{1}{2}, \qquad
    - \frac{\tau}{2} \leqq \Im \lam_j \leqq \frac{\tau}{2}, \qquad
\end{equation}
without loss of generality thanks to the following lemma.

\begin{lem}
\label{Bv-period}
Suppose $(\nu, \{\lam_1, \dots, \lam_M\})$ is a solution of the Bethe
equations. Then for any $j$, $1\leqq j \leqq M$,

i) $\Psi_\nu(\lam_1,\dots, \lam_j + 1, \dots, \lam_M)$ is proportional
to $\Psi_\nu(\lam_1, \dots, \lam_j, \dots, \lam_M)$, and
$(\nu, \{\lam_1, \dots, \lam_j + 1, \dots, \lam_M\})$ is
a solution of the Bethe equations.

ii) $\Psi_{\nu+2}(\lam_1,\dots, \lam_j + \tau, \dots, \lam_M)$ is
proportional to $\Psi_\nu(\lam_1, \dots, \lam_j, \dots, \lam_M)$, and
$(\nu+2, \{\lam_1, \dots, \lam_j + \tau, \dots, \lam_M\})$ is
a solution of the Bethe equations.
\end{lem}

\begin{pf}
The quasi-periodicity of theta functions implies that
$(\nu,   \{\lam_1, \dots, \lam_j + 1, \dots, \lam_M\})$ and
$(\nu+2, \{\lam_1, \dots, \lam_j + \tau, \dots, \lam_M\})$ are solutions
of the Bethe equations.

Substituting $\lam \to \lam +1$ and $\lam \to \lam + \tau$ in
\eqref{def:L-n}, we obtain
\begin{align*}
    L_n(\lam + 1   ) &= - \sig^1 L_n(\lam) \sig^1,
\\
    L_n(\lam + \tau) &= - e^{- \pi i \tau - 2 \pi i \lam}
                        \sig^3 L_n(\lam) \sig^3.
\end{align*}
In the same way (see \eqref{def:M-k}):
\begin{align*}
    M_k(\lam + 1   ;s) &= \sig^1 M_k(\lam;s),
\\
    M_k(\lam + \tau;s) &= \sig^3 M_k(\lam;s)
                          \begin{pmatrix}
                          -e^{ \pi i (s_+ +2k\eta-\lam-\tau)} & 0\\
                          0 & -e^{-\pi i (s_- +2k\eta+\lam)} \\
                          \end{pmatrix}.
\end{align*}
Hence from \eqref{def:twist-monod} follows
\begin{gather}
    \calT_{k,k'}(\lam+1   ;s) = (-1)^N \calT_{k,k'}(\lam),
\label{tw-monod:+1}
\\
\begin{split}
    \calT_{k,k'}(\lam+\tau;s) =& (-1)^N e^{\pi iN \tau - 2\pi i N\lam}
                          \begin{pmatrix}
                          e^{- \pi i (s_+ + 2k\eta -\lam-\tau)} & 0\\
                          0 & e^{\pi i (s_- + 2k\eta +\lam)}
                          \end{pmatrix} \times\\
                          \times &
                          \calT_{k,k'} (\lam;s)
                          \begin{pmatrix}
                          e^{ \pi i (s_+ + 2k'\eta -\lam-\tau)} & 0\\
                          0 & e^{-\pi i (s_- + 2k'\eta +\lam)}
                          \end{pmatrix}.
\end{split}
\label{tw-monod:+tau}
\end{gather}
The (1,2)-component of \eqref{tw-monod:+1} is
$B_{k,k'}(\lam+1) = (-1)^N B_{k,k'}(\lam)$. Thus
\begin{equation*}
    \Psi_\nu(\lam_1, \dots, \lam_j + 1, \dots, \lam_M) =
    (-1)^N \Psi_\nu(\lam_1, \dots, \lam_j, \dots, \lam_M),
\end{equation*}
which proves i).
The (1,2)-component of \eqref{tw-monod:+tau} is
\begin{equation*}
    B_{k,k'}(\lam+\tau) = \text{ const. } e^{-2 \pi i (k+k') \eta}
    B_{k,k'}(\lam).
\end{equation*}
Here const.~does not depend on $k$, $k'$, but depends on $\lam$,
$s_\pm$. Thus
\begin{equation*}
    \Phi_a(\lam_1, \dots, \lam_j+\tau, \dots, \lam_M) =
    \text{ const. } e^{-4 \pi i a \eta}
    \Phi_a(\lam_1, \dots, \lam_j, \dots, \lam_M),
\end{equation*}
and
\begin{equation*}
\begin{split}
    \Psi_{\nu+2}(\lam_1, \dots, \lam_j+\tau, \dots, \lam_M) &=
    \text{ const. } \sum_{a=0}^{r-1}
    e^{2 \pi i(\nu+2)a \eta - 4 \pi ia \eta}
    \Phi_a(\lam_1, \dots, \lam_M)
\\
    &= \text{ const. } \Psi_\nu(\lam_1, \lam_2, \dots, \lam_M).
\end{split}
\end{equation*}
This proves ii).
\end{pf}

Baxter developed the coordinate Bethe Ansatz in \cite{bax:eight},
expanding eigenvectors into linear combination of fundamental vectors
$|a_0, \dots, a_N \rangle$. We have such a coordinate expression
of the above defined Bethe vectors. Namely
\begin{prop}
Let $\Psi_\nu(\lam_1, \dots, \lam_M)$ be as defined by \eqref{Bethe-vec}.
Then
\begin{equation}
\begin{split}
    \Psi_\nu(\lam_1, \dots, \lam_M)
    =&
    \sum_{a=0}^{r-1} e^{2 \pi i \nu \eta a} \times \\
    \times
    \sum\begin{Sb}a_0 = a, a_1, \dots \\ \dots, a_{N-1}, a_N = a\end{Sb}
    &\left(
    \sum_{ \{a_{i,j}\} }
    \prod_{i=1}^M \prod_{j=1}^N
    \faceW{a_{i,j}}{a_{i,j-1}}{a_{i-1,j}}{a_{i-1,j-1}}{\lam_i}
    \right)
    |a_N, \dots, a_1, a_0 \rangle.
\end{split}
\label{coord-expr}
\end{equation}
Here the sum in $()$ is taken over a set of integers $a_{i,j}$
($0\leqq i \leqq M$, $0 \leqq j \leqq N$) satisfying
{\em the admissibility condition}
\begin{equation*}
    a_{i,j} - a_{i-1,j} = \pm 1, \qquad
    a_{i,j} - a_{i,j-1} \in
    \{ -2\ell, -2\ell + 2, \dots, 2\ell - 2, 2\ell\},
\end{equation*}
and {\em the boundary condition},
\begin{gather*}
    a_{0, j} = a_j,  \qquad a_{M, j} = a - N\ell + 2\ell j,
\\
    a_{i, 0} = a - i,\qquad a_{i, N} = a + i.
\end{gather*}
\end{prop}

Note that the Bethe equations are not assumed here.

\begin{pf}
Operator $B_{a_N,a_0}(\lam)$ is the $(1,2)$-element of the monodromy
matrix $\calT_{a_N,a_0}(\lam) = L_{a_N,a_{N-1}} \dots L_{a_1,a_0}$.
Hence by \eqref{vertex-face'}:
\begin{align*}
    B_{a_N,a_0}(\lam) & |a_N, a_{N-1}, \dots, a_1, a_0 \rangle
    =\\
    =&
    \sum\begin{Sb}\eps_N = -, \eps_{N-1}, \dots \\
                  \dots, \eps_1, \eps_0 = +\end{Sb}
    L_{a_N, a_{N-1}}(\eps_N, \eps_{N-1};\lam) \dots
    L_{a_j, a_{j-1}}(\eps_j, \eps_{j-1};\lam) \dots \\
    \dots &
    L_{a_1, a_0    }(\eps_1, \eps_0    ;\lam) \cdot
    \phi_{a_N,a_{N-1}} \tensor \dots \tensor
    \phi_{a_j,a_{j-1}} \tensor \dots \tensor
    \phi_{a_1,a_0}
\\
    =&
    \sum\begin{Sb}\eps_N = -, \eps_{N-1}, \dots \\
                  \dots, \eps_1, \eps_0 = +\end{Sb}
    \prod_{j=1}^N
    \faceW{a_j}{a_{j-1}}{a_j+\eps_j}{a_{j-1}+\eps_{j-1}}{\lam} \cdot \\
    &\cdot
    \phi_{a_N+\eps_N,a_{N-1}+\eps_{N-1}} \tensor \dots \tensor
    \phi_{a_j+\eps_j,a_{j-1}+\eps_{j-1}} \tensor \dots \tensor
    \phi_{a_1+\eps_1,a_0    +\eps_0    }.
\end{align*}
Applying this formula iteratively, we arrive at \eqref{coord-expr}.
\end{pf}

\subsection{Sum rule}
\setcounter{equation}{0}
\label{sumrule}
In the previous section we defined Bethe vectors by \eqref{Bethe-vec}.
Here we show an integrality condition of sum of parameters
$\lam_j$.

\begin{thm}
\label{sumrule:thm}
Let $(\nu, \{\lam_1, \dots, \lam_M\})$ be a solution of the Bethe
equations \eqref{Bethe-eq}. We assume that $\{\lam_j\}$ satisfy the
following additional conditions: For any $j= 1, \dots, M$,

i) $\lam_j \not\in \{ 2(n+\ell)\eta \,|\, n \in \Integer\}$.

ii) there exists $a \in \Integer$ such that
$\lam_j + 2 a\eta \not\equiv \lam_k$ $\pmod{\Integer+\Integer\tau}$
for any $k=1, \dots, M$.

iii) (Technical assumption of non-degeneracy: see
\appref{sumrule:pf}.)

Then there exist integers $n_0$, $n_1$ which satisfy
\begin{equation}
    2 \sum_{j=1}^M \lam_j = n_0 + n_1 \tau.
\label{sumrule:eq}
\end{equation}
\end{thm}

The proof is technical and contained in \appref{sumrule:pf}.
Note that assumptions i) and ii) in \thmref{sumrule:thm} are satisfied
under the string hypothesis in \secref{string}. Assumption iii) is
hard to check. For string solutions considered in \secref{gr.state}
and \secref{s-matrix}, \eqref{sumrule:eq} is checked directly.

Baxter derived this rule in \cite{bax:part}, directly constructing
an operator on $\calH$ which gives $Q(\lam)$ (see \eqref{def:Q}) as its
eigenvalue. Unfortunately we have not yet found such an operator in our
context. (Kulish and Reshetikhin \cite{kul-resh} found for a rational $R$
matrix case that transfer matrix ``converges'' to the $Q$ operator by
iteration of fusion procedures.)

In addition, Baxter's result tells us that
$\sum_{j=1}^M \lam_j$ is related to parities of Bethe vectors. These
parities are associated to reversing the arrows ($\sigma^1$) or to
assigning $-1$ to the down arrows ($\sigma^3$) of the $XYZ$ spin chain
model. They correspond to the unitary operators $U_1$ and $U_3$ acting
on the spin $\ell$ representations. (See \appref{skl-alg}.)

\begin{lem}
For $a=1,2,3$, $U_a$ commutes with the $L$ operator as
\begin{equation}
    U_a^{-1} L(\lam) U_a = (\sig^a)^{-1} L(\lam) \sig^a,
\label{comm-rel:U-L}
\end{equation}
They commute with the transfer matrix:
$U_a^{-1} T(\lam) U_a = T(\lam)$.
\end{lem}

\begin{pf}
Adjoint action by $U_a$ induces an automorphism $X_a$ on the Sklyanin
algebra (see \appref{skl-alg}):
\begin{align*}
    U_a^{-1} L(\lam) U_a &= X_a(L(\lam))
\\
    &= W_0(\lam) \rho^\ell (S^0) \tensor \sig^0
     + W_a(\lam) \rho^\ell (S^a) \tensor \sig^a \\
    &- W_b(\lam) \rho^\ell (S^b) \tensor \sig^b
     - W_c(\lam) \rho^\ell (S^c) \tensor \sig^c,
\end{align*}
where $(a,b,c)$ is a cyclic permutation of $(1,2,3)$.
By the anti-commutativity of the Pauli matrices, the right hand side
of the above equation is nothing but $(\sig^a)^{-1} L(\lam) \sig^a$.

Commutativity $T(\lam) U_a = U_a T(\lam)$ is a direct consequence
of \eqref{comm-rel:U-L} and \eqref{def:transfer}.
\end{pf}

Operators $U_a^{\tensor N}$ on $\calH$ are involutive and commute with
each other. In fact, by virtue of relations $U_a^2 = (-1)^{2\ell}$ and
$U_a U_b = (-1)^{2\ell} U_b U_a$, we have
\begin{gather*}
    (U_a^{\tensor N})^2 = (-1)^{2N\ell} = 1,
\\
                  U_a^{\tensor N} U_b^{\tensor N} =
    (-1)^{2N\ell} U_b^{\tensor N} U_a^{\tensor N} =
                  U_b^{\tensor N} U_a^{\tensor N}.
\end{gather*}
(Recall that $N\ell$ is an integer.) Therefore an eigenvalue of
$U_a^{\tensor N}$ is either $+1$ or $-1$. Assume that $\Psi \in \calH$
is a common eigenvector of $T(\lam)$ and $U_a$'s. We assign parities
$\nu''$ and $\nu'$ to $\Psi$ by
\begin{equation*}
    U_1^{\tensor N} \Psi = (-1)^{\nu''} \Psi, \qquad
    U_3^{\tensor N} \Psi = (-1)^{\nu'}  \Psi.
\end{equation*}

{}From Baxter's result \cite{bax:part} and \thmref{sumrule:thm},
the following conjecture seems to be plausible.

\begin{conj}
Let $(\nu, \{ \lam_1, \dots, \lam_M \})$ be a solution of the
Bethe equations, and $\nu''$ and $\nu'$ be parities of the Bethe
vector $\Psi_\nu(\lam_1, \dots,\lam_M)$ defined above. Then
\begin{align*}
                 \nu + \nu' + N\ell &\equiv 0 \quad \pmod{2},
\\
    \nu\tau - 2 \sum_{j=1}^M \lam_j &\equiv \nu'' + N\ell
    \quad\pmod{2}.
\end{align*}
\end{conj}

%
%
%
%
\section{Thermodynamic limit}

In this chapter we consider the limit $N \to \infty$.
Our calculation is based on the string hypothesis introduced by
Takahashi and Suzuki \cite{tak-suz} which we review in
\secref{string}. In \secref{gr.state} we compute the free energy of the
model. In \secref{s-matrix} we introduce four kinds of perturbation of
the string configuration of the ground state. Each of them are
parametrized by two continuous parameters which are regarded as
rapidities of physical particles. We compute polarization of the Dirac
sea of quasi-particles induced by these perturbation, following the
recipe by Johnson, Krinsky and McCoy \cite{jo-kr-mc}. We also
calculate  eigenvalues of the $S$ matrix of two physical particles,
using the method developed by Korepin \cite{kor}, Destri and
Lowenstein \cite{des-low}. The result coincides with Smirnov's
conjecture \cite{fijkmy}.

\subsection{String hypothesis}
\setcounter{equation}{0}
\label{string}
First let us rescale the parameters so that integrals in later
sections are taken over a segment in the real line. We denote
$x_j = i t \lam_j$. Then the Bethe equations \eqref{Bethe-eq} takes
the form
\begin{equation}
    \left(\frac{\theta_{11}(x_j + 2i\ell\eta t ;it)}
               {\theta_{11}(x_j - 2i\ell\eta t ;it)}\right)^N
    =
    e^{-4\pi i \eta (\nu + 2 \sum_{k=1}^M x_k)}
    \prod\begin{Sb}k=1 \\ k \neq j \end{Sb}^M
    \frac{\theta_{11}(x_j - x_k + 2i\eta t ;it)}
         {\theta_{11}(x_j - x_k - 2i\eta t ;it)},
\label{Bethe-eq'}
\end{equation}
while the corresponding eigenvalue is
\begin{equation}
\begin{aligned}
    t(x) &= (-2i\sqrt{t})^N
    e^{\frac{\pi N}{t}(x^2 - 4 \ell(\ell+1) \eta^2 t^2)} \times \\
    &\left(
    e^{ 2\pi i \eta (\nu + 2 \sum_{j=1}^M x_j)}
    (2\theta_{11}(x + 2i\ell\eta t ;it))^N
    \prod_{j=1}^M
    \frac{\theta_{11}(x - x_j - 2i\eta t ;it)}
         {\theta_{11}(x - x_j;            it)} \right.\\
    &\left. +
    e^{-2\pi i \eta (\nu + 2 \sum_{j=1}^M x_j)}
    (2\theta_{11}(x - 2i\ell\eta t ;it))^N
    \prod_{j=1}^M
    \frac{\theta_{11}(x - x_j + 2i\eta t ;it)}
         {\theta_{11}(x - x_j;            it)} \right),
\end{aligned}
\label{Bethe-ev'}
\end{equation}
where $x=it\lam$. According to \lemref{Bv-period}, we may assume that
$|\Re(x_j)| \leqq 1/2$ for any $j$.

Now {\em the string hypothesis} \cite{tak-suz} is stated in the following
way. For sufficiently large $N$ solutions of \eqref{Bethe-eq'} cluster
into groups known as $A$-strings, $A = 1, 2, \dots$, with parity $\pm$:
\begin{equation}
    x^{A,\pm}_{j,\alpha}= x^{A,\pm}_j + 2 i\eta t \alpha
                        + O(e^{-\delta N}),
    \qquad
    \alpha = \frac{-A+1}{2}, \frac{-A+3}{2}, \dots, \frac{A-1}{2},
\label{def:A-string}
\end{equation}
where $\Im x^{A,+}_j \equiv 0$ $\pmod{\Integer t}$ and
$\Im x^{A,-}_j \equiv t/2$ $\pmod{\Integer t}$. Complex numbers
$x^{A,\pm}_j$ is called a center of a string. Due to \lemref{Bv-period}
we may assume that  $\Im x^{A,+}_j = 0$ and $\Im x^{A,-}_j = t/2$. We
denote the number of $A$-strings with parity $\pm$ by $\sharp(A,\pm)$.

Note that assumption i) of \thmref{sumrule:thm} is satisfied for
$A$-strings with parity $+$, if $A \equiv 2\ell$ $\pmod{2}$ and for
any strings with non-zero real abscissa. Assumption ii) holds for
$A$-strings, $A < r$, if real parts of centers of all strings lie in
the interval $[-1/2, 1/2]$.

Since we are interested only in the thermodynamic limit $N\to\infty$,
we neglect the exponentially small deviation term $O(e^{-\delta N})$
in what follows.

\subsection{Ground state and free energy}
\setcounter{equation}{0}
\label{gr.state}
The result in this section was announced in \cite{take2}.
The ground state configuration is specified as follows: $\nu=0$,
$\sharp(2\ell, +) = N/2$, $\sharp(A,\pm) = \sharp(2\ell, -) =0$
for $A\neq 2\ell$ and centers of $2\ell$-strings distribute
symmetrically around 0. (Hence $\sum_{j=1}^{N/2} x_j^{2\ell,+} = 0$).
This is consistent with the result of the $XXX$
type model by Takhtajan \cite{takh}, Babujian \cite{bab}, that of the
$XXZ$ type by Sogo \cite{sog}, Kirillov and Reshetikhin
\cite{kir-resh} and that of the $XYZ$ model by Baxter \cite{bax:part}.

Multiplying the Bethe equations \eqref{Bethe-eq'} for
$x_j = x^{2\ell,+}_{j,\alpha}$,
$\alpha = -\ell + 1/2, \dots, \ell - 1/2$ (cf.~\eqref{def:A-string}),
and taking the logarithm, we obtain
\begin{multline}
    N \sum_{\alpha=-\ell + 1/2}^{\ell - 1/2}
      \Phi(x_j; 2i(\alpha + \ell)\eta t) =\\
    =
    2\pi Q^{2\ell}_j +
    \sum_{k=1}^{N/2} \left(
    \sum_{m=1}^{2\ell-1} \Phi(x_j-x_k; 2im\eta t) +
    \sum_{m=0}^{2\ell-1} \Phi(x_j-x_k; 2i(m+1)\eta t)
    \right),
\label{log-Be:gr.st}
\end{multline}
where we omit the index $(2\ell,+)$ of $x^{2\ell,+}_j$ and function
$\Phi$ is defined by:
\begin{equation}
    \Phi(x;i\mu t) = \frac{1}{i}
                   \log \frac{\theta_{11}(x + i\mu t ;it)}
                             {\theta_{11}(x - i\mu t ;it)} + \pi.
\label{def:Phi}
\end{equation}
We take the branch of $\Phi$ so that $\Phi(\pm 1/2;i\mu) = \mp\pi$,
$\Phi(0;i\mu) =0$. Half integers $Q^{2\ell}_j$ in \eqref{log-Be:gr.st}
specify the branches of logarithm. Applying Takhtajan-Faddeev's
argument \cite{takh-fad:xxx} to our case, they satisfy
$Q^{2\ell}_j - Q^{2\ell}_{j-1} = -1$.
We also assume that $x_j$ are ordered by $j$: $x_j > x_{j-1}$.
(Note that $\Phi(x;i\mu t)$ is a decreasing function by
\lemref{positivity}.)

We assume that centers of $2\ell$ strings fill in the interval
$(-1/2, 1/2)$ with density $\rho(x)$ in the limit $N\to\infty$,
\begin{equation*}
    \frac{1}{N(x_{j+1} - x_j)} \to \rho(x), \quad x_j\to x,
    \qquad N\to\infty.
\end{equation*}
Subtracting \eqref{log-Be:gr.st} for $j$ from that for $j+1$ and taking
the limit, we obtain an integral equation
\begin{multline}
    \sum_{\alpha=-\ell + 1/2}^{\ell - 1/2}
    \Phi'(x; 2i(\alpha + \ell)\eta t) = - 2\pi \rho(x) + \\
    +
    \int_{-1/2}^{1/2} \left(
    \sum_{m=1}^{2\ell-1} \Phi'(x-y; 2im\eta t) +
    \sum_{m=0}^{2\ell-1} \Phi'(x-y; 2i(m+1)\eta t)
    \right) \rho(y) \,dy.
\label{int-eq:rho}
\end{multline}
This equation is easily solved by Fourier expansion. Using formula
\eqref{Fourier:Phi'}, we have
\begin{equation}
    \rho(x) = \sum_{n\in\Integer}
              \frac{e^{2\pi i nx}}{2\ch 2 \pi n\eta t}.
\label{rho}
\end{equation}

Let us compute the eigenvalue for this Bethe vector.
The expression of the eigenvalue \eqref{Bethe-ev'} consists of two terms:
\begin{equation}
\begin{aligned}
    \Lam_1(x) =& (-2i\sqrt{t})^N
    e^{\frac{\pi N}{t}(x^2 - 4 \ell(\ell+1) \eta^2 t^2)} \times\\
    \times &
    (2\theta_{11}(x + 2i\ell\eta t ;it))^N
    \prod_{j=1}^M
    \frac{\theta_{11}(x - x_j - 2i\eta t ;it)}
         {\theta_{11}(x - x_j;            it)},
\\
    \Lam_2(x) =& (-2i\sqrt{t})^N
    e^{\frac{\pi N}{t}(x^2 - 4 \ell(\ell+1) \eta^2 t^2)} \times\\
    \times &
    (2\theta_{11}(x - 2i\ell\eta t ;it))^N
    \prod_{j=1}^M
    \frac{\theta_{11}(x - x_j + 2i\eta t ;it)}
         {\theta_{11}(x - x_j;            it)},
\end{aligned}
\end{equation}
and the eigenvalue is $t(x) = \Lam_1(x) + \Lam_2(x)$.
Both $\Lam_1$ and $\Lam_2$ contribute equally to $t(x)$ in
the thermodynamic limit $N\to\infty$ for $\ell > 1/2$, since
\begin{equation*}
\begin{split}
    \frac{1}{N}\log \frac{\Lam_1}{\Lam_2} &=
    \log \frac{\theta_{11}(x + 2i \ell\eta t ;it)}
              {\theta_{11}(x - 2i \ell\eta t :it)} +\\
    &+ \frac{1}{N}
    \sum_{j=1}^{N/2} \left(
    \frac{\theta_{11}(x - x_j - i(2\ell+1)\eta t ;it)}
         {\theta_{11}(x - x_j + i(2\ell+1)\eta t ;it)}
    +
    \frac{\theta_{11}(x - x_j - i(2\ell-1)\eta t ;it)}
         {\theta_{11}(x - x_j + i(2\ell-1)\eta t ;it)}
    \right)
\\
    \overset{N\to\infty}{\longrightarrow}&
    \log \frac{\theta_{11}(x + 2i \ell\eta t ;it)}
              {\theta_{11}(x - 2i \ell\eta t :it)} +\\
    +&
    \int_{-1/2}^{1/2}
    \left(
    \frac{\theta_{11}(x - y - i(2\ell+1)\eta t ;it)}
         {\theta_{11}(x - y + i(2\ell+1)\eta t ;it)}
    +
    \frac{\theta_{11}(x - y - i(2\ell-1)\eta t ;it)}
         {\theta_{11}(x - y + i(2\ell-1)\eta t ;it)}
    \right) \times \\
    & \qquad\qquad \times \rho(y) \, dy
\\
    =& -2\pi i.
\end{split}
\end{equation*}
In the case of the eight vertex model ($\ell = 1/2$) for $\lam>0$,
$\Lam_1$ is dominant in magnitude, because
\begin{equation}
\begin{split}
    \frac{1}{N}\log \frac{\Lam_1}{\Lam_2} &=
    \log \frac{\theta_{11}(x + i \eta t ;it)}
              {\theta_{11}(x - i \eta t :it)}
    + \frac{1}{N}
    \sum_{j=1}^{N/2}
    \frac{\theta_{11}(x - x_j - 2i\eta t ;it)}
         {\theta_{11}(x - x_j + 2i\eta t ;it)}
\\
    &\overset{N\to\infty}{\longrightarrow}
    - \frac{3\pi i}{2} - \pi i x
    - \sum_{n=1}^\infty
      \frac{i \sin 2\pi n x}{n \ch 2 \pi n \eta t},
\end{split}
\end{equation}
the real part of which is positive. (Recall that $x=it\lam$,
$\lam>0$.) This is a subtle difference between $\ell = 1/2$ and higher
spin cases, but the final result does not differ much.
Namely, in the thermodynamic limit,
\begin{multline}
    \frac{1}{N}\log t(x)
    \overset{N\to\infty}{\longrightarrow}
    \frac{1}{N}\log \Lam_1(x) \\
    \overset{N\to\infty}{\sim}
    \log(-2i) + \frac{1}{2}\log t +
    \frac{\pi}{t}(x^2 - 4\ell(\ell+1)\eta^2 t^2) +
    \log 2 +\\
    +
    \log \theta_{11}(x+2i\ell\eta t ;it) +
    \int_{-1/2}^{1/2}
    \log \frac{\theta_{11}(x - y - 2i\eta t ;it)}
              {\theta_{11}(x - y;            it)} \rho(y) \, dy.
\end{multline}
Substituting \eqref{rho} into this, we obtain the free energy
\eqref{free-energy}:
\begin{multline}
    - \beta f(\lam) =
    \text{(const.) } +
    \log \theta_{11}( \lam+2\ell\eta ;\tau)
    - 2 \pi t (\lam - \eta)(1 - 4\ell \eta)\\
    -\sum_{n=1}^\infty
    \frac
    {\sh \pi nt(1 - 4\ell\eta)\,\sh 2\pi nt(\lam-\eta)}
    {n \sh \pi nt\, \ch 2\pi n\eta t}.
\end{multline}
Here (const.) is an unessential term which does not depend on $\lam$.

\subsection{Low-lying excitations and S matrices}
\setcounter{equation}{0}
\label{s-matrix}

As is seen in \secref{gr.state}, the ground state consists of $N/2$
$2\ell$-strings filling the Dirac sea. In this section we perturb this
Dirac sea, slightly changing the string configuration. Since we are
interested in the two particle states, we choose such configurations
that reduce to two particle states of models associated to
trigonometric and rational $R$ matrices \cite{takh-fad:xxx},
\cite{takh}, \cite{kir-resh}, \cite{bab}, \cite{sog}.

Let us consider the following configurations:
\begin{enumerate}
\renewcommand{\labelenumi}{(\Roman{enumi})}
\item $\sharp(2\ell  , +) = N/2-2$,
      $\sharp(2\ell-1, +) = 1$,
      $\sharp(2\ell+1, +) = 1$;
\item $\sharp(2\ell  , +) = N/2-1$,
      $\sharp(2\ell-1, +) = 1$,
      $\sharp(1      , -) = 1$.
\end{enumerate}
We call the Bethe vectors specified by these data {\em excited state} I
and II respectively. In the higher spin $XXX$ case \cite{takh}, for
example, two particle states are specified by similar configurations;
one (singlet) is the same as I above, the other (triplet) is
defined by $\sharp(2\ell,+) = N/2-1$, $\sharp(2\ell-1,+)=1$. Since
one-string with parity $-$ goes away to infinity when $t$ tends to
$\infty$, we can expect that excited state II reduces to the triplet
state in the rational limit. As a matter of course, when $\ell = 1/2$,
($2\ell-1$)-string is absent. Hence the following argument needs to be
modified, but one obtains results for $\ell=1/2$ by simply
putting $\ell=1/2$ in formulae for general $\ell$. We do not mention
this modification.

\subsubsection*{Excited state I}
Now we consider the excited state I. We omit the plus sign designating
the parity, since all string have parity $+$.
Multiplying the Bethe equations \eqref{Bethe-eq'} for a $2\ell$-string
$x_j = x^{2\ell}_{j,\alpha}$, $\alpha = -\ell + 1/2, \dots, \ell - 1/2$
with the center $x^{2\ell}_j$, and taking the logarithm, we obtain
\begin{equation}
\begin{split}
    N \sum_{\alpha=-\ell + 1/2}^{\ell - 1/2} &
      \Phi(x^{2\ell}_j; 2i(\alpha + \ell)\eta t)
    =
    2\pi Q^{2\ell}_j +
    8\pi \ell \eta (\nu + 2\Sigma_{\text I}) +\\
    +&
    \sum_{k=1}^{N/2-2} \left(
    \sum_{m=1}^{2\ell-1} \Phi(x^{2\ell}_j-x^{2\ell}_k; 2im\eta t) +
    \sum_{m=0}^{2\ell-1} \Phi(x^{2\ell}_j-x^{2\ell}_k; 2i(m+1)\eta t)
    \right)\\
    +&
    \sum_{m=1/2}^{2\ell - 3/2} (
    \Phi(x^{2\ell}_j-x_-^{2\ell-1}; 2im\eta t) +
    \Phi(x^{2\ell}_j-x_-^{2\ell-1}; 2i(m+1)\eta t)
    )\\
    +&
    \sum_{m=1/2}^{2\ell - 1/2} (
    \Phi(x^{2\ell}_j-x_+^{2\ell+1}; 2im\eta t) +
    \Phi(x^{2\ell}_j-x_+^{2\ell+1}; 2i(m+1)\eta t)
    ),
\end{split}
\label{log-Be:I-2l}
\end{equation}
where $x_\pm^{2\ell\pm 1}$ are the centers of the
$(2\ell\pm1)$-strings,
\begin{equation}
    \Sigma_{\text I} = \text{ (sum of all } x^A_{j,\alpha}) =
    2\ell \sum_{j=1}^{N/2-2} x^{2\ell}_j +
    (2\ell-1) x_-^{2\ell-1} +
    (2\ell+1) x_+^{2\ell+1}.
\label{sum:I}
\end{equation}
The argument of \cite{takh-fad:xxx} applied to \eqref{log-Be:I-2l}
implies that there are $N/2$ vacancies for $Q^{2\ell}_j$'s. Thus there remain
two vacancies (holes) left unoccupied by centers of $2\ell$-strings.

We renumber the centers of strings as follows.
\begin{enumerate}
\renewcommand{\labelenumi}{(\roman{enumi})}
\item $2\ell$-strings: $x_j$, $j=1, \dots, N/2$, $j\neq j_1, j_2$,
where $Q^{2\ell}_{j_1}$ and $Q^{2\ell}_{j_2}$ correspond to holes. Following
the argument in \cite{takh-fad:xxx} again, we assume that
$x_j > x_{j'}$ if $j>j'$.
\item $2\ell-1$-string: $x_-$.
\item $2\ell+1$-string: $x_+$.
\end{enumerate}

In the thermodynamic limit centers of $2\ell$-strings fill
the interval $(-1/2, 1/2)$ continuously with density $\rho_{\text I}(x)$
and two holes at $x_1 = \lim x_{j_1}$ and $x_2 = \lim x_{j_2}$.
(We abuse indices.) Subtracting \eqref{log-Be:I-2l} for $j$ from
that for $j+1$, we obtain
\begin{equation}
\begin{split}
    \sum_{\alpha=-\ell + 1/2}^{\ell - 1/2} &
    \Phi'(x; 2i(\alpha + \ell)\eta t) =
    - 2\pi \left(
    \rho_{\text I}(x) + \frac{1}{N} (\delta(x-x_1) + \delta(x-x_2))
    \right) +\\
    +&
    \int_{-1/2}^{1/2} \left(
    \sum_{m=1}^{2\ell-1} \Phi'(x-y; 2im\eta t) +
    \sum_{m=0}^{2\ell-1} \Phi'(x-y; 2i(m+1)\eta t)
    \right) \rho_{\text I}(y)\, dy\\
    +& \frac{1}{N}
    \sum_{m=1/2}^{2\ell - 3/2}
    (\Phi'(x-x_-; 2im\eta t) + \Phi'(x-x_-; 2i(m+1)\eta t))\\
    +& \frac{1}{N}
    \sum_{m=1/2}^{2\ell - 1/2}
    (\Phi'(x-x_+; 2im\eta t) + \Phi'(x-x_+; 2i(m+1)\eta t)),
\end{split}
\label{int-eq:I-2l}
\end{equation}
for large $N$.
The solution of this integral equation for $\rho_{\text I}(x)$ is
\begin{equation}
    \rho_{\text I}(x) = \rho(x)
              + \frac{1}{N}(\sig(x-x_1) + \sig(x-x_2)
                          + \omega_-(x-x_-) + \omega_+(x-x_+)),
\label{def:rho-I}
\end{equation}
where $\rho(x)$ is defined above, $\sig(x)$ and $\omega_\pm(x)$ are
solutions of the following integral equations:
Integral equation for $\sig(x)$:
\begin{multline}
    2\pi \sig(x) = -2\pi \delta(x) \\
    +
    \int_{-1/2}^{1/2} \left(
    \sum_{m=1}^{2\ell-1} \Phi'(x-y; 2im\eta t) +
    \sum_{m=0}^{2\ell-1} \Phi'(x-y; 2i(m+1)\eta t)
    \right) \sig(y) \, dy,
\label{int-eq:sigma}
\end{multline}
Integral equation for $\omega_-(x)$:
\begin{multline}
    2\pi \omega_-(x) = \\
    \int_{-1/2}^{1/2} \left(
    \sum_{m=1}^{2\ell-1} \Phi'(x-y; 2im\eta t) +
    \sum_{m=0}^{2\ell-1} \Phi'(x-y; 2i(m+1)\eta t)
    \right) \omega_-(y) \, dy \\
    +
    \sum_{m=1/2}^{2\ell - 3/2}
    (\Phi'(x; 2im\eta t) + \Phi'(x; 2i(m+1)\eta t)),
\label{int-eq:omega-}
\end{multline}
Integral equation for $\omega_+(x)$:
\begin{multline}
    2\pi \omega_+(x) = \\
    \int_{-1/2}^{1/2} \left(
    \sum_{m=1}^{2\ell-1} \Phi'(x-y; 2im\eta t) +
    \sum_{m=0}^{2\ell-1} \Phi'(x-y; 2i(m+1)\eta t)
    \right) \omega_+(y) \, dy \\
    +
    \sum_{m=1/2}^{2\ell - 1/2}
    (\Phi'(x; 2im\eta t) + \Phi'(x; 2i(m+1)\eta t)).
\label{int-eq:omega+}
\end{multline}
They are easily solved by the Fourier expansion explicitly:
\begin{align}
    \sig(x) &= -\frac{1}{4\ell}
    - \sum_{n=1}^\infty
    \frac{\sh \pi nt\,\sh 2 \pi n\eta t}
         {\sh \pi nt(1-4\ell\eta)\,\sh 4 \pi n\ell\eta t\,
          \ch 2\pi n\eta t}
    \cos 2 \pi nx,
\label{def:sigma}
\\
    \omega_-(x) &= -\frac{2\ell-1}{2\ell}
    - \sum_{n=1}^\infty
    \frac{2 \sh 2 \pi n(2\ell-1)\eta t}{\sh 4 \pi n\ell\eta t}
    \cos 2 \pi nx,
\label{def:omega-}
\\
    \omega_+(x) &= -1
    - \sum_{n=1}^\infty
    \frac{2 \sh \pi nt (1-2(2\ell+1)\eta)}{\sh \pi nt(1-4\ell\eta)}
    \cos 2 \pi nx.
\label{def:omega+}
\end{align}

Product of the Bethe equations \eqref{Bethe-eq'} for
the ($2\ell-1$)-string $x_- + 2 i \alpha \eta t$,
$\alpha = -\ell+1, \dots, \ell-1$, gives the equation:
\begin{equation}
\begin{split}
    N \sum_{\alpha = -\ell+1}^{\ell+1}&
    \Phi(x_-; 2i (\alpha+\ell)\eta t) =
    2\pi Q^{2\ell-1}_- - (2\ell-1) 4 \pi \eta (\nu + 2\Sigma_{\text I})
    \\
    +&
    \sum_{k=1, k\neq j_1, j_2}^{N/2}
    \sum_{m=1/2}^{2\ell - 3/2}
    (\Phi(x_- - x_k; 2im\eta t) + \Phi(x_- - x_k; 2i(m+1)\eta t))\\
    +&
    \sum_{m=1}^{2\ell - 1}
    (\Phi(x_- - x_+; 2im\eta t) + \Phi(x_- - x_+; 2i(m+1)\eta t))
\end{split}
\label{log-Be:I-2l-1}
\end{equation}
This time there exists only one vacancy for $Q^{2\ell-1}_-$ which
determines the branch. We set $Q^{2\ell-1}_- = 0$.
In the thermodynamic limit equation \eqref{log-Be:I-2l-1} gives
an integral equation:
\begin{equation}
\begin{split}
    \sum_{\alpha = -\ell+1}^{\ell+1}&
    \Phi(x_-; 2i (\alpha+\ell)\eta t) =
    - \frac{1}{N} (2\ell-1) 4 \pi \eta (\nu + 2\Sigma_{\text I})\\
    +&
    \int_{-1/2}^{1/2}
    \sum_{m=1/2}^{2\ell - 3/2}
    (\Phi(x_- - y; 2im\eta t) + \Phi(x_- - y; 2i(m+1)\eta t))
    \rho_{\text I}(y)\, dy\\
    +& \frac{1}{N}
    \sum_{m=1}^{2\ell - 1}
    (\Phi(x_- - x_+; 2im\eta t) + \Phi(x_- - x_+; 2i(m+1)\eta t)).
\end{split}
\label{int-eq:I-2l-1}
\end{equation}
This equation reduces to
\begin{equation}
    \frac{2\ell}{2\ell-1}
    \int^{x_- - x_1}_{- x_- + x_2} \omega_-(y)\,dy
    + x_+ - \frac{x_1 + x_2}{2}
    = (1-4\ell\eta) \Sigma_{\text I} - 2\ell \nu \eta
\label{sumrule:I-2l-1}
\end{equation}
by \eqref{def:rho-I}, \eqref{rho}, \eqref{def:sigma},
\eqref{def:omega-}, \eqref{def:omega+} and \eqref{Fourier:Phi}.

On the other hand, product of the Bethe equations \eqref{Bethe-eq'}
for the ($2\ell+1$)-string $x_+ + 2 i \alpha \eta t$,
$\alpha = -\ell, \dots, \ell$, gives the equation:
\begin{equation}
\begin{split}
    N \sum_{\alpha = -\ell+1}^{\ell}&
    \Phi(x_+; 2i (\alpha+\ell)\eta t) =
    2\pi Q^{2\ell+1}_+ - (2\ell+1) 4 \pi \eta (\nu + 2\Sigma_{\text I})
    \\
    +&
    \sum_{k=1, k \neq j_1, j_2}^{N/2}
    \sum_{m=1/2}^{2\ell - 1/2}
    (\Phi(x_+ - x_k; 2im\eta t) + \Phi(x_+ - x_k; 2i(m+1)\eta t))\\
    +&
    \sum_{m=1}^{2\ell - 1}
    (\Phi(x_+ - x_-; 2im\eta t) + \Phi(x_+ - x_-; 2i(m+1)\eta t))
\end{split}
\label{log-Be:I-2l+1}
\end{equation}
Again only one vacancy for $Q^{2\ell+1}_+$ which determines the branch
exists. We set $Q^{2\ell+1}_+ = 0$. The integral equation in the
thermodynamic limit given by \eqref{log-Be:I-2l+1} is
\begin{equation}
\begin{split}
    \sum_{\alpha = -\ell+1}^{\ell}&
    \Phi(x_+; 2i (\alpha+\ell)\eta t) =
    - \frac{1}{N} (2\ell+1) 4 \pi \eta (\nu + 2\Sigma_{\text I})\\
    +&
    \int_{-1/2}^{1/2}
    \sum_{m=1/2}^{2\ell - 1/2}
    (\Phi(x_+ - y; 2im\eta t) + \Phi(x_+ - y; 2i(m+1)\eta t))
    \rho_{\text I}(y)\, dy\\
    +& \frac{1}{N}
    \sum_{m=1}^{2\ell - 1}
    (\Phi(x_+ - x_-; 2im\eta t) + \Phi(x_+ - x_-; 2i(m+1)\eta t)).
\end{split}
\label{int-eq:I-2l+1}
\end{equation}
and hence
\begin{equation}
    \frac{1}{2}
    \int^{x_+ - x_1}_{- x_+ + x_2} \omega_+(y)\,dy
    + x_+ - \frac{x_1 + x_2}{2} =
    (1-2(2\ell+1)\eta) \Sigma_{\text I} - (2\ell+1) \nu \eta
\label{sumrule:I-2l+1}
\end{equation}
by \eqref{def:rho-I}, \eqref{rho}, \eqref{def:sigma}, \eqref{def:omega-},
\eqref{def:omega+} and \eqref{Fourier:Phi}.

Let us denote the solution of the Bethe equations for the ground state
by $x_j^{G}$. {\em Polarization} of the Dirac sea of $2\ell$-strings
for excited state I is defined by
\begin{equation}
    J(x) := \rho(x) \lim_{N\to\infty} N(x_j - x_j^G)
          = \lim_{N\to\infty} \frac{x_j-x_j^G}{x_{j+1}^G-x_j^G},
\label{def:pol}
\end{equation}
where $x=\lim_{N\to\infty} x_j$. (See \cite{kor}, \cite{jo-kr-mc}.)
Subtracting \eqref{log-Be:I-2l} from \eqref{log-Be:gr.st} and using
the integral equation \eqref{int-eq:rho}, one can derive the integral
equation for $J(x)$:
\begin{equation}
\begin{split}
    -2\pi J(x)
    +&
    \int_{-1/2}^{1/2} \left(
    \sum_{m=1}^{2\ell-1} \Phi'(x-y; 2im\eta t) +
    \sum_{m=0}^{2\ell-1} \Phi'(x-y; 2i(m+1)\eta t)
    \right) J(y)\, dy\\
    =&
    - 8\pi \ell \eta( \nu + 2 \Sigma_{\text I} ) \\
    +&
    \sum_{m=1/2}^{2\ell - 3/2}
    (\Phi(x-x_-; 2im\eta t) + \Phi(x-x_-; 2i(m+1)\eta t))\\
    +&
    \sum_{m=1/2}^{2\ell - 1/2}
    (\Phi(x-x_+; 2im\eta t) + \Phi(x-x_+; 2i(m+1)\eta t))\\
    -&
    \sum_{a=1,2} \left(
    \sum_{m=1}^{2\ell-1} \Phi'(x-x_a; 2im\eta t) +
    \sum_{m=0}^{2\ell-1} \Phi'(x-x_a; 2i(m+1)\eta t)
    \right).
\end{split}
\label{int-eq:I-pol}
\end{equation}
Thus the polarization is determined as
\begin{gather}
    J(x) = \sum_{n\in\Integer} J_n e^{2\pi i n x},
\\
    J_0 = \eta (\nu + 2 \Sigma_{\text I})
        - \frac{2\ell-1}{2\ell} x_- - x_+
        + \frac{4\ell-1}{4\ell} (x_1 + x_2),
\label{I:J-0:pol}
\\
\begin{aligned}
    J_n &= \frac{\sh 2\pi n(2\ell-1)\eta t}
               {2\pi i n\sh 4\pi n\ell\eta t}
          \left(e^{-2\pi in x_-} -
          \frac{e^{-2\pi in x_1} + e^{-2\pi in x_2}}
               {2\ch 2\pi n\eta t}
          \right)\\
         &+
          \frac{\sh \pi n(1-2(2\ell+1)\eta) t}
               {2\pi i n\sh \pi n(1-4\ell\eta)t}
          \left(e^{-2\pi in x_+} -
          \frac{e^{-2\pi in x_1} + e^{-2\pi in x_2}}
               {2\ch 2\pi n\eta t}
          \right).
\end{aligned}
\end{gather}
On the other hand, by the definition of the polarization \eqref{def:pol},
\begin{equation}
    2\ell \int_{-1/2}^{1/2} J(x)\,dx
    =
    \Sigma_{\text I} - (2\ell+1) x_+ - (2\ell-1) x_- + 2\ell(x_1 + x_2).
\label{I:J-0:sum}
\end{equation}
Combining \eqref{I:J-0:pol} and \eqref{I:J-0:sum}, we obtain
\begin{equation}
    x_+ - \frac{x_1 + x_2}{2} =
    (1-4\ell\eta) \Sigma_{\text I} - 2\ell\nu\eta.
\label{sumrule:I-pol}
\end{equation}

Now we determine $x_\pm$ in terms of $x_1$, $x_2$ regarded as free
parameters. We have derived three equations connecting $x_1$, $x_2$
and $x_\pm$: \eqref{sumrule:I-2l-1}, \eqref{sumrule:I-2l+1} and
\eqref{sumrule:I-pol}.
{}From \eqref{sumrule:I-2l-1} and \eqref{sumrule:I-pol} follows
\begin{equation}
    \int^{x_- - x_1}_{- x_- + x_2} \omega_-(y)\,dy
    = 0.
\label{sumrule:x-}
\end{equation}
Thus $x_-= (x_1+x_2)/2$, since $\omega_-(y) < 0$ because of
\eqref{def:omega-} and \lemref{positivity}.
Equations \eqref{sumrule:I-2l+1} and \eqref{sumrule:I-pol} imply
\begin{equation}
    \frac{1}{2}
    \int^{x_+ - x_1}_{- x_+ + x_2} \omega_+(y)\,dy
    = -\eta(\nu + 2\Sigma_{\text I}),
\label{x+:1}
\end{equation}
or, equivalently,
\begin{equation}
    \int^{x_+ - x_1}_{- x_+ + x_2}
    \left(\omega_+(y) + \frac{2\eta}{1-4\ell\eta}\right)\,dy
    = -\nu \frac{2\eta}{1-4\ell\eta}.
\label{x+:2}
\end{equation}
Equation \eqref{x+:2} shows that $x_+$ is uniquely determined for
each $\nu$ because of \lemref{positivity}. Constraints
\eqref{sumrule:I-pol}, \eqref{x+:2} and \eqref{sumrule:eq} on $\nu$,
$x_+$ and $\Sigma_{\text I}$ are simultaneously satisfied if we put
\begin{equation*}
    \nu = k \left( \frac{r}{2} - (2\ell+1)r' \right), \qquad
    \Sigma_{\text I} = (2\ell+1)\frac{k r'}{2}, \qquad
    x_+ = \frac{x_1 + x_2}{2} + \frac{k r'}{2},
\end{equation*}
where $k$ is an arbitrary integer. Recall that $r$ is assumed to be
even. Apparently there are infinitely many solutions, but in fact
only two Bethe vectors are independent in this series:
\begin{lem}
For two integers $k,k'$ such that $k\equiv k'$ $\pmod{2}$,
corresponding Bethe vectors are linearly dependent.
\end{lem}
\begin{pf}
Let us denote $\nu$ and $x_+$ corresponding to $k$ and $k'$ by
($\nu(k),x_+(k)$) and ($\nu(k'),x_+(k')$) respectively. Then
\begin{align*}
    \nu(k) - \nu(k') &= - (2\ell+1) (k-k')r' + \frac{k-k'}{2} r\\
                     &\equiv - (2\ell+1) (k-k')r' \qquad\pmod{r},
\\
    x_+(k) - x_+(k') &= \frac{k-k'}{2} r'.
\end{align*}
Shifting $x_+$ by one means shifting $2\ell +1$ of solutions of Bethe
equations by one. Lemma is proved by \lemref{Bv-period}.
\end{pf}

If we take $k=0$, then $\nu=0$, $\Sigma_{\text I}=0$ and $x_+ = (x_1+x_2)/2$.
We call the Bethe vector corresponding to this configuration
{\em excited state} I${}_0$.

Since $r$ and $r'$ are coprime, there exists an integer $k$ such that
$kr'\equiv 1$ $\pmod{r}$. Shifting $x_+$ by an integer, we may assume
that $x_+ = (x_1 + x_2)/2 + 1/2$ with a suitable $\nu$. (See
\lemref{Bv-period}.) We call this Bethe vector {\em excited state}
I${}_1$.

It seems that there are no other solutions for $x_+$ and $\nu$, since
the integrality conditions ($\nu\in\Integer$ and \thmref{sumrule:thm})
are very strong.

\subsubsection*{Excited state II}
Now we consider the excited state II. As in the case of the excited
state I, multiplying the Bethe equations \eqref{Bethe-eq'} for a
$2\ell$-string $x_i = x^{2\ell,+}_{j,\alpha}$,
$\alpha = -\ell + 1/2, \dots, \ell - 1/2$ with the center
$x^{2\ell,+}_j$, and taking the logarithm, we obtain
\begin{equation}
\begin{split}
    N &\sum_{\alpha=-\ell + 1/2}^{\ell - 1/2}
      \Phi(x^{2\ell,+}_j; 2i(\alpha + \ell)\eta t) =
    2\pi Q^{2\ell}_j +
    8\pi \ell \eta (\nu + 2\Sigma_{\text{II}}) \\
    +&
    \sum_{k=1}^{N/2-1} \left(
    \sum_{m=1}^{2\ell-1} \Phi(x^{2\ell,+}_j-x^{2\ell,+}_k; 2im\eta t) +
    \sum_{m=0}^{2\ell-1} \Phi(x^{2\ell,+}_j-x^{2\ell,+}_k; 2i(m+1)\eta t)
    \right)\\
    +&
    \sum_{m=1/2}^{2\ell - 3/2} (
    \Phi(x^{2\ell,+}_j-x^{2\ell-1,+}_-; 2im\eta t) +
    \Phi(x^{2\ell,+}_j-x^{2\ell-1,+}_-; 2i(m+1)\eta t)
    )\\
    +&
    \Psi(x^{2\ell,+}_j-x_0; 2i(2\ell+1)\eta t) +
    \Psi(x^{2\ell,+}_j-x_0; 2i(2\ell-1)\eta t)
    ),
\end{split}
\label{log-Be:II-2l}
\end{equation}
where $x^{2\ell-1,+}_-$ is the center of the $(2\ell-1)$-string,
$\{x_0 + it/2\}$ ($x_0\in\Real$) is the one-string with parity $-$,
\begin{equation}
    \Sigma_{\text{II}}
    = \text{ (sum of all } x^A_{j,\alpha}) - \frac{it}{2}
    =
    2\ell \sum_{j=1}^{N/2-1} x^{2\ell,+}_j +
    (2\ell-1) x^{2\ell-1,+} +
    x_0^{1,-},
\label{sum:II}
\end{equation}
and function $\Psi(x;i\mu t)$ is defined by
\begin{equation}
    \Psi(x;i\mu t) = \frac{1}{i}
                   \log \frac{\theta_{01}(x + i\mu t ;it)}
                             {\theta_{01}(x - i\mu t ;it)}.
\label{def:Psi}
\end{equation}

By the same argument as for the excited state I,
there are $N/2+1$ vacancies for $Q^{2\ell}_j$'s. Thus there are again
two holes of centers of $2\ell$-strings.

We renumber the centers of strings as follows.
\begin{enumerate}
\renewcommand{\labelenumi}{(\roman{enumi})}
\item $2\ell$-strings: $x_j$, $j=1, \dots, N/2+1$, $j\neq j_1, j_2$,
where $Q^{2\ell}_{j_1}$ and $Q^{2\ell}_{j_2}$ correspond to
holes. Following the argument in \cite{takh-fad:xxx} again, we assume
that $x_j > x_{j'}$ if $j>j'$. The string with its center at
$x_{N/2+1}$ will be placed at the zone boundary $x=1/2$ in the
thermodynamic limit.
\item ($2\ell-1$)-string: $x_-$.
\item $1$-string with parity $-$: $x_0 + \frac{it}{2}$.
\end{enumerate}

As in the previous case, we obtain an integral equation for the
density of centers of $2\ell$-strings, $\rho_{\text{II}}(x)$, on the
interval $(-1/2, 1/2)$
\begin{equation}
\begin{split}
    \sum_{\alpha=-\ell + 1/2}^{\ell - 1/2} &
    \Phi'(x^{2\ell}_j; 2i(\alpha + \ell)\eta t) =
    - 2\pi \left(
    \rho_{\text{II}}(x) + \frac{1}{N} (\delta(x-x_1) + \delta(x-x_2))
    \right) \\
    +&
    \int_{-1/2}^{1/2} \left(
    \sum_{m=1}^{2\ell-1} \Phi'(x-y; 2im\eta t) +
    \sum_{m=0}^{2\ell-1} \Phi'(x-y; 2i(m+1)\eta t)
    \right) \rho_{\text{II}}(y)\, dy\\
    +& \frac{1}{N}
    \sum_{m=1/2}^{2\ell - 3/2}
    (\Phi'(x-x_-; 2im\eta t) + \Phi'(x-x_-; 2i(m+1)\eta t))\\
    +& \frac{1}{N}
    \Psi'(x-x_0; 2i(2\ell+1)\eta t) + \Psi'(x-x_0; 2i(2\ell-1)\eta t)
\end{split}
\label{int-eq:II-2l}
\end{equation}
for large $N$.
Its solution is
\begin{equation}
    \rho_{\text{II}}(x) = \rho(x)
              + \frac{1}{N}(\sig(x-x_1) + \sig(x-x_2)
              + \omega_-(x-x_-) + \omega_0(x-x_0)),
\label{def:rho-II}
\end{equation}
where $\rho(x)$ (\eqref{int-eq:rho}, \eqref{rho}), $\sig(x)$
(\eqref{int-eq:sigma}, \eqref{def:sigma}), $\omega_-(x)$
(\eqref{int-eq:omega-}, \eqref{def:omega-}) are as defined before,
and $\omega_0(x)$ is a solution of the following integral equation:
\begin{equation}
\begin{split}
    2\pi \omega_0(x) =&
    \int_{-1/2}^{1/2} \left(
    \sum_{m=1}^{2\ell-1} \Phi'(x-y; 2im\eta t) +
    \sum_{m=0}^{2\ell-1} \Phi'(x-y; 2i(m+1)\eta t)
    \right) \omega_0(y) \, dy \\
    +&
    \Psi'(x; i(2\ell+1)\eta t) + \Psi'(x; i(2\ell-1)\eta t).
\label{int-eq:omega0}
\end{split}
\end{equation}
Explicitly, $\omega_0(x)$ is
\begin{equation}
    \omega_0(x) =
    \sum_{n=1}^\infty
    \frac{2 \sh 2 \pi n\eta t}{\sh \pi n(1-4\ell\eta)t}
    \cos 2 \pi nx.
\label{def:omega0}
\end{equation}

The Bethe equations \eqref{Bethe-eq'} for the ($2\ell-1$)-string gives
the equation:
\begin{equation}
\begin{split}
    N \sum_{\alpha = -\ell+1}^{\ell+1}&
    \Phi(x_-; 2i (\alpha+\ell)\eta t) =
    2\pi Q^{2\ell-1}_- - (2\ell-1) 4 \pi \eta (\nu + 2\Sigma_{\text{II}})
    \\
    +&
    \sum_{k=1, \neq j_1, j_2}^{N/2+1}
    \sum_{m=1/2}^{2\ell - 3/2}
    (\Phi(x_- - x_k; 2im\eta t) + \Phi(x_- - x_k; 2i(m+1)\eta t))\\
    +&
    \Psi(x_- - x_0; 2i\ell\eta t) + \Psi(x_- - x_0; 2i(\ell - 1)\eta t)).
\end{split}
\label{log-Be:II-2l-1}
\end{equation}
We set $Q^{2\ell-1}_- = 0$, since there is only one vacancy.
The corresponding integral equation in the thermodynamic limit is:
\begin{equation}
\begin{split}
    \sum_{\alpha = -\ell+1}^{\ell+1} &
    \Phi(x_-; 2i (\alpha+\ell)\eta t) =
    - \frac{1}{N} (2\ell-1) 4 \pi \eta (\nu + 2\Sigma_{\text{II}})\\
    +&
    \int_{-1/2}^{1/2}
    \sum_{m=1/2}^{2\ell - 3/2}
    (\Phi(x_- - y; 2im\eta t) + \Phi(x_- - y; 2i(m+1)\eta t))
    \rho_{\text{II}}(y)\, dy\\
    +& \frac{1}{N}
    (\Psi(x_- - x_0; 2i\ell\eta t) + \Psi(x_- - x_0; 2i(\ell-1)\eta t)).
\end{split}
\label{int-eq:II-2l-1}
\end{equation}
This equation reduces to
\begin{multline}
    \frac{\ell}{2\ell-1}
    \int^{x_- - x_1}_{- x_- + x_2} \omega_-(y)\,dy
    + x_0 - \frac{x_1 + x_2}{2}
    = (1-4\ell\eta) \Sigma_{\text{II}} - 2 \ell \nu \eta.
\label{sumrule:II-2l-1}
\end{multline}
as before.

The Bethe equations \eqref{Bethe-eq'}
for the 1-string with parity $-$ gives the equation:
\begin{equation}
\begin{split}
    N \Psi(x_0; 2i \ell\eta t) =&
    2\pi Q^{1-}_0 - 4 \pi \eta (\nu + 2\Sigma_{\text{II}})\\
    +&
    \sum_{k=1, \neq j_1, j_2}^{N/2}
    (\Psi(x_0 - x_k; i(2\ell+1)\eta t) +
     \Psi(x_0 - x_k; i(2\ell-1)\eta t))\\
    +&
    \Psi(x_0 - x_-; 2i\ell\eta t) + \Psi(x_0 - x_-; 2i(\ell-1)\eta t)
\end{split}
\label{log-Be:II-0}
\end{equation}
Again we choose the branch $Q^{1-}_0 = 0$. The integral equation is
\begin{equation}
\begin{split}
    \Psi(x_0; 2i \ell\eta t) =&
    - \frac{1}{N} 4 \pi \eta (\nu + 2\Sigma_{\text{II}})\\
    +&
    \int_{-1/2}^{1/2}
    (\Psi(x_0 - y; i(2\ell+1)\eta t) + \Psi(x_0 - y; i(2\ell-1)\eta t))
    \rho_{\text{II}}(y)\, dy\\
    +& \frac{1}{N}
    (\Psi(x_0 - x_-; 2i\ell\eta t) + \Psi(x_0 - x_-; 2i(\ell-1)\eta t)).
\end{split}
\label{int-eq:II-0}
\end{equation}
This is rewritten as
\begin{equation}
    \int^{x_0 - x_1}_{- x_0 + x_2} \omega_0(y)\,dy
    = -2\eta( 2\Sigma_{\text{II}} + \nu).
\label{sumrule:II-0}
\end{equation}

Let us compute the polarization $J(x)$ of this state.
Subtracting \eqref{log-Be:II-2l} from \eqref{log-Be:gr.st} and using
the integral equation \eqref{int-eq:rho}, we can derive the integral
equation for $J(x)$:
\begin{equation}
\begin{split}
    -2\pi J(x)
    +&
    \int_{-1/2}^{1/2} \left(
    \sum_{m=1}^{2\ell-1} \Phi'(x-y; 2im\eta t) +
    \sum_{m=0}^{2\ell-1} \Phi'(x-y; 2i(m+1)\eta t)
    \right) J(y)\, dy\\
    =&
    - 8\pi \ell \eta( \nu + 2 \Sigma_{\text{II}} )\\
    +&
    \sum_{m=1/2}^{2\ell - 3/2}
    (\Phi(x-x_-; 2im\eta t) + \Phi(x-x_-; 2i(m+1)\eta t))\\
    +&
    \Psi(x-x_0; i(2\ell+1)\eta t) + \Psi(x-x_0; i(2\ell-1)\eta t)\\
    +&
    \sum_{m=1}^{2\ell-1} \Phi'(x-\half; 2im\eta t) +
    \sum_{m=0}^{2\ell-1} \Phi'(x-\half; 2i(m+1)\eta t)\\
    -&
    \sum_{a=1,2} \left(
    \sum_{m=1}^{2\ell-1} \Phi'(x-x_a; 2im\eta t) +
    \sum_{m=0}^{2\ell-1} \Phi'(x-x_a; 2i(m+1)\eta t)
    \right).
\end{split}
\label{int-eq:II-pol}
\end{equation}
Hence the polarization is
\begin{gather}
    J(x) = \sum_{n\in\Integer} J_n e^{2\pi i n x},
\\
    J_0 = -\frac{x}{4\ell}+ \eta (\nu + 2 \Sigma_{\text{II}}) - \frac{1}{2}
        - \frac{2\ell-1}{2\ell} x_-
        + \frac{4\ell-1}{4\ell} (x_1 + x_2),
\label{II:J-0:pol}
\\
\begin{aligned}
    J_n &= \frac{\sh 2\pi n(2\ell-1)\eta t}{\sh 4\pi n\ell\eta t}
        \times\\
        &\times
          \left(e^{-2\pi in x_-} -
          \frac{- e^{-\pi i n} + e^{-2\pi in x_1} + e^{-2\pi in x_2}}
               {2\ch 2\pi n\eta t}
          -
          \frac{e^{-2\pi i n x_0}}{\sh \pi n t(1-4\ell\eta)}
          \right)\\
         &+
          \frac{\sh \pi n(1-2(2\ell+1)\eta) t}
               {2\ch 2\pi n\eta t \sh \pi n(1-4\ell\eta)t}
          (- e^{-\pi i n} + e^{-2\pi in x_1} + e^{-2\pi in x_2})
\end{aligned}
\end{gather}
On the other hand, by the definition of the polarization \eqref{def:pol},
\begin{equation}
    2\ell \int_{-1/2}^{1/2} J(x)\,dx
    =
    \Sigma_{\text{II}} - (2\ell-1) x_- - x_0 - \ell + 2\ell(x_1 + x_2).
\label{II:J-0:sum}
\end{equation}
It follows from \eqref{II:J-0:sum} and \eqref{II:J-0:pol} that
\begin{equation}
    x_0 - \frac{x_1 + x_2}{2} =
    (1-4\ell\eta) \Sigma_{\text{II}} - 2\ell\nu\eta.
\label{sumrule:II-pol}
\end{equation}

{}From \eqref{sumrule:II-2l-1} and \eqref{sumrule:II-pol} follows
the same equation as \eqref{sumrule:x-} and thus $x_-= (x_1+x_2)/2$ as in
the case of the excited state I.
Equations \eqref{sumrule:II-0} and \eqref{sumrule:II-pol} imply
\begin{equation}
    \int^{x_0 - x_1}_{- x_0 + x_2}
    \left(\omega_0(y) + \frac{2\eta}{1-4\ell\eta}\right)\,dy
    = -\nu \frac{2\eta}{1-4\ell\eta}.
\label{x0}
\end{equation}
This equation is uniquely solved because of \eqref{def:omega0} and
\lemref{positivity}: $x_0 = (x_1 + x_2)/2 - \nu/2$. \lemref{Bv-period}
tells that only two cases give independ Bethe vectors:
($\nu=0$, $x_0 = (x_1+x_2)/2$) which we call the excited state
II${}_0$ and ($\nu=1$, $x_0=(x_1+x_2)/2 + \half$) which we call
the excited state II${}_1$.

\begin{table}[h]
\caption{Two particle excited states}
\begin{center}
\begin{tabular}{c|l|l|l|l}
\hline
Excited & $2\ell$-strings & $2\ell+1$-string & $2\ell-1$-string &
1-string \\
States  & parity $+$      & parity $+$       & parity $+$       &
parity $-$ \\
\hline
I${}_0$ & density $\rho_{\text I}$ & $x_+=$        & $x_-=$        &\\
        & holes $x_1$, $x_2$       & $(x_1+x_2)/2$ & $(x_1+x_2)/2$ &\\
\hline
I${}_1$ & density $\rho_{\text I}$ & $x_+=$          & $x_-=$          &\\
        & holes $x_1$, $x_2$       & $(x_1+x_2+1)/2$ & $(x_1+x_2)/2$ &\\
\hline
II${}_0$ & density $\rho_{\text{II}}$ && $x_-=$        & $x_0=$       \\
         & holes $x_1$, $x_2$         && $(x_1+x_2)/2$ & $(x_1+x_2)/2$\\
\hline
II${}_1$ & density $\rho_{\text{II}}$ && $x_-=$        & $x_0=$         \\
         & holes $x_1$, $x_2$         && $(x_1+x_2)/2$ & $(x_1+x_2+1)/2$\\
\hline
\end{tabular}
\end{center}
\end{table}

\subsubsection*{$S$ matrix}

Above we found four excited states with two free parameters
$x_1$, $x_2$: I${}_0$, I${}_1$, II${}_0$, II${}_1$.
In the rational limit, $t \to \infty$, $\eta\to 0$, $\eta t$ fixed,
the string configuration of I${}_0$ becomes that of the singlet state
of the corresponding spin chain, whereas the configurations I${}_1$,
II${}_0$, II${}_1$ seem to approach to that of the triplet states,
since the one-string with parity $-$ and the string with abscissa
$x_+=(x_1+x_2+1)/2$ goes beyond the sight. (Recall that real abscissas
of strings are rescaled so that they fill the whole real line in the
limit.) Hence one might expect that these four Bethe vectors give four
dimensional space of two physical particle states (spin waves) of the
corresponding spin chain. In fact for the eight vertex model
\begin{equation}
\begin{split}
    &
    \log T(x)|_{\text{excited state}} -
    \log T(x)|_{\text{ground state}} \\
    =&
    \log \tau(x - i\eta t - x_1) +
    \log \tau(x - i\eta t - x_2),
\end{split}
\end{equation}
where (excited state) means any one of the excited states
I${}_0$, I${}_1$, II${}_0$, II${}_1$ and
\begin{equation}
    \log\tau(x) := -\frac{\pi i}{2} - \pi i x
                   - i \sum_{n=1}^\infty
                   \frac{\sin 2\pi n x}{n \ch 2\pi n \eta t},
\end{equation}
(see \cite{jo-kr-mc} for details of calculations).
This means that all conserved quantities such as momentum $P(x)$
or energy over the ground states are split into two terms:
$$
    P(x) = -\pi x - \sum_{n=1}^\infty
                   \frac{\sin 2\pi n x}{n \ch 2\pi n \eta t},
$$
and thus we can regard these excited states as two particle states of
the $XYZ$ spin chain.

For higher spin cases we have not yet computed fused transfer matrix
which corresponds to the spin chain with local interaction.
{}From the result of the rational and trigonomeric models
\cite{takh}, \cite{kir-resh}, we conjecture that the momentum and the
energy of physical particles do not depend on the spin $\ell$.
Based on this conjecture, we calculate the $S$ matrix of two physical
particles.

%

As is discussed in \cite{dmn}, the $S$ matrix of physical particles (spin
waves) could depend on the way of calculation for the case of higher spin.
In order to make our standpoint clear, let us recall the calculation of
eigenvalues of the $S$ matrix in more details, following \cite{kor} and
Section 9 of \cite{low} (cf.~also \cite{des-low}, \cite{takh-fad:xxx}):
The appearance of an integer $Q_j^{2\ell}$ in
\eqref{log-Be:I-2l} and \eqref{log-Be:II-2l}
can be interpreted as a consequence of the periodic boundary condition
which we imposed on the lattice. Namely, if we move a physical particle
around the whole chain, the total phase shift of the wave function
accumlated should be an integer multiple of $2\pi i$. The main
contribution comes from the momentum $P$ of the particle as
$i P N$, the {\it free phase\/}:
$$
    i P N =
    -iN \pi x
    -iN \sum_{n=1}^\infty
    \frac{\sin 2\pi n x}{n \ch 2\pi n \eta t}.
$$
Note that the right hand side of the above equation is eventually expressed
as the logarithm of (the right hand side)$/$(the left hand side) of the
Bethe equation for the ground state \eqref{log-Be:gr.st}:
\begin{equation}
\begin{split}
    &- iN \sum_{\alpha=-\ell + 1/2}^{\ell - 1/2}
      \Phi(x; 2i(\alpha + \ell)\eta t) \\
    &+   i \sum_{k=1}^{N/2} \left(
    \sum_{m=1}^{2\ell-1} \Phi(x-x_k; 2im\eta t) +
    \sum_{m=0}^{2\ell-1} \Phi(x-x_k; 2i(m+1)\eta t)
    \right)
\\
    &\overset{N\to\infty}{\longrightarrow}
    - iN \sum_{\alpha=-\ell + 1/2}^{\ell - 1/2}
      \Phi(x; 2i(\alpha + \ell)\eta t) \\
    &+
    iN \int_{-1/2}^{1/2} \left(
    \sum_{m=1}^{2\ell-1} \Phi(x-y; 2im\eta t) +
    \sum_{m=0}^{2\ell-1} \Phi(x-y; 2i(m+1)\eta t)
    \right) \rho(y) \,dy
\\
    &= iPN.
\end{split}
\label{free-phase}
\end{equation}
Hence the ground state can be interpreted as a Dirac sea of non-interacting
particles, since the momenta of particles are integer multiples of
$2\pi/N$ because of \eqref{log-Be:gr.st} and the above equation.

For the excited state, however, the phase shift comes not only
from this free phase but also from the interaction between physical
particles. Because of the periodic boundary condition which fixes the total
phase shift to an integer multiple of $2\pi i$, this means that the
calculation of scattering phase shift of a physical particle is equivalent
to the calculation of $O(1/N)$ shift of the momentum. In other words, the
$S$ matrix of physical particles can be calculated by splitting the total
phase shift, an integer multiple of $2\pi i$, into the free phase $i P N$
of order $O(N)$ and the scattering phase of order $O(1)$.

We consider the excited state I first. The total phase shift for the
physical particle with rapidity $x_1$ can be read off from
\eqref{log-Be:I-2l} as follows:
\begin{equation}
\begin{split}
    -N &\sum_{\alpha=-\ell + 1/2}^{\ell - 1/2}
      \Phi(x_1; 2i(\alpha + \ell)\eta t)
    +
    8\pi \ell \eta (\nu + 2\Sigma_{\text I}) +\\
    +&
    \sum_{k=1}^{N/2-2} \left(
    \sum_{m=1}^{2\ell-1} \Phi(x_1-x^{2\ell}_k; 2im\eta t) +
    \sum_{m=0}^{2\ell-1} \Phi(x_1-x^{2\ell}_k; 2i(m+1)\eta t)
    \right)\\
    +&
    \sum_{m=1/2}^{2\ell - 3/2} (
    \Phi(x_1-x_-^{2\ell-1}; 2im\eta t) +
    \Phi(x_1-x_-^{2\ell-1}; 2i(m+1)\eta t)
    )\\
    +&
    \sum_{m=1/2}^{2\ell - 1/2} (
    \Phi(x_1-x_+^{2\ell+1}; 2im\eta t) +
    \Phi(x_1-x_+^{2\ell+1}; 2i(m+1)\eta t)
    ),
\end{split}
\label{total-phase:I}
\end{equation}
which is equal to an integer multiple of $2\pi i$ because of
\eqref{log-Be:I-2l}. Subtracting the free phase contribution $iPN$
\eqref{free-phase} from the total phase \eqref{total-phase:I}, and taking
the limit $N\to\infty$, we obtain the remainder of order $O(1)$ and thus we
can interpret it as the scattering phase shift from the above argument.
An explicit expression for the eigenvalue of the $S$ matrix for
the excited state I is as follows:
\begin{equation}
\begin{split}
    i \log &(\pm S_{\text I}(x_1 - x_2)) =
    -8\pi \ell \eta (\nu + 2\Sigma_{\text I})\\
    +&
    N \int_{-1/2}^{1/2} \left(
    \sum_{m=1}^{2\ell-1} \Phi(x_1-y; 2im\eta t) +
    \sum_{m=0}^{2\ell-1} \Phi(x_1-y; 2i(m+1)\eta t)
    \right)\times \\
    &\qquad \times (\rho_{\text I}(y)-\rho(y)) \, dy\\
    +&
    \sum_{m=1/2}^{2\ell - 3/2}
    (\Phi(x_1-x_-; 2im\eta t) + \Phi(x_1-x_-; 2i(m+1)\eta t))\\
    +&
    \sum_{m=1/2}^{2\ell - 1/2}
    (\Phi(x_1-x_+; 2im\eta t) + \Phi(x_1-x_+; 2i(m+1)\eta t)),
\end{split}
\label{def:log-S:I}
\end{equation}
The term $-8\pi \ell \eta \nu$ can be interpreted as an effect
from the background or boundary, while the rest of the right hand side comes
from interaction of pseudo-particles. The ambiguity of sign comes from
normalizations of asymptotic states. The right hand side is computed by
integrating \eqref{int-eq:sigma}, \eqref{int-eq:omega-},
\eqref{int-eq:omega+}. The result is:
\begin{equation}
\begin{split}
    i \log &(\pm S_{\text I}(x)) =\\
    =&\sum_{n=1}^\infty
    \left(
    \frac{\sh \pi n t (1-4\ell\eta - 2\eta)}
         {n \sh \pi n t (1-4\ell\eta)\, \ch 2\pi n \eta t} +
    \frac{\sh \pi n t (4\ell\eta - 2\eta)}
         {n \sh 4 \pi n \ell\eta t \, \ch 2\pi n \eta t}
    \right)
    \sin 2\pi n x \\
    +&
    \sum_{n=1}^\infty
    \frac{2 \sh \pi nt (4\ell \eta - 2\eta)}
         {\sh 4 \pi n \ell \eta t}
    \sin \pi n x \\
    +&
    \sum_{n=1}^\infty
    \frac{\sh \pi n t(2\eta -(1-4\ell\eta))}
         {\sh \pi n t(1-4\ell\eta)}
    \sin \pi n (x - \eps),
\end{split}
\label{log-S:I}
\end{equation}
where $\eps$ is 0 or 1 for the excited state I${}_0$ or I${}_1$,
respectively. The first term of the right hand side come from
holes ($\sig(x-x_2)$ in $\rho_{\text I}(x) - \rho(x)$ of
\eqref{def:log-S:I}), the second term from the ($2\ell-1$)-string
($\omega_-(x-x_-)$) and the last term from the ($2\ell+1$)-string
($\omega_+(x-x_+)$).

Computation for the excited state II is the same.
The result is
\begin{equation}
\begin{split}
    i \log &(\pm S_{\text{II}}(x)) =\\
    = &\sum_{n=1}^\infty
    \left(
    \frac{\sh \pi n t (1-4\ell\eta - 2\eta)}
         {n \sh \pi n t (1-4\ell\eta)\, \ch 2\pi n \eta t}
    +
    \frac{\sh \pi n t (4\ell\eta - 2\eta)}
         {n \sh 4 \pi n \ell\eta t \, \ch 2\pi n \eta t}
    \right) \sin 2\pi n x \\
    +&
    \sum_{n=1}^\infty
    \frac{2 \sh \pi nt (4\ell\eta - 2\eta)}
         {\sh 4 \pi n \ell \eta t}
    \sin \pi n x \\
    +& \pi + \pi x
    +
    \sum_{n=1}^\infty
    \frac{\sh 2 \pi n \eta t}
         {\sh \pi n t(1-4\ell\eta)}
    \sin \pi n (x - \eps),
\end{split}
\label{log-S:II}
\end{equation}
where $\eps$ is 0 or 1 for the excited state II${}_0$ or II${}_1$,
respectively. As in the case of the excited state I, the first term
of the right hand side come from holes, the second term from the
($2\ell-1$)-string and the last term from the one-string with parity $-$
($\omega_0(x-x_0)$).

We fix the signs left undetermined so that the above $S$ matrix is
the permutation matrix in the non-interacting limit $x=0$ and the
excited states I${}_0$ reduces to a singlet while other three states
form a triplet as in the rational and trigonometric cases.
Then,
\begin{equation}
\begin{aligned}
    S(x)|_{\text{excited state I${}_0$}}
    &= S_0(x)
    \frac{\theta_{11}\left(\frac{x}{2}-it\eta;it(1-4\ell\eta)\right)}
         {\theta_{11}\left(\frac{x}{2}+it\eta;it(1-4\ell\eta)\right)},
\\
    S(x)|_{\text{excited state I${}_1$}}
    &= S_0(x)
    \frac{\theta_{10}\left(\frac{x}{2}-it\eta;it(1-4\ell\eta)\right)}
         {\theta_{10}\left(\frac{x}{2}+it\eta;it(1-4\ell\eta)\right)},
\\
    S(x)|_{\text{excited state II${}_0$}}
    &= S_0(x)
    \frac{\theta_{01}\left(\frac{x}{2}-it\eta;it(1-4\ell\eta)\right)}
         {\theta_{01}\left(\frac{x}{2}+it\eta;it(1-4\ell\eta)\right)},
\\
    S(x)|_{\text{excited state II${}_1$}}
    &= S_0(x)
    \frac{\theta_{00}\left(\frac{x}{2}-it\eta;it(1-4\ell\eta)\right)}
         {\theta_{00}\left(\frac{x}{2}+it\eta;it(1-4\ell\eta)\right)},
\end{aligned}
\label{spec(S)}
\end{equation}
where
\begin{equation}
    S_0(x) = e^{-2\pi i x}
    \frac{\theta_{11}\left(\frac{x}{2}-it\eta;4\ell it\eta\right)}
         {\theta_{11}\left(\frac{x}{2}+it\eta;4\ell it\eta\right)}
    \bbbS(x; 1-4\ell\eta) \bbbS(x; 4\ell\eta).
\end{equation}
and function $\bbbS(x;\mu)$ is defined by
\begin{equation}
\begin{split}
    \bbbS(i\lam t; \mu) &=
    \exp \left(
    \sum_{n=1}^\infty
    \frac{\sh \pi n t(\mu - 2\eta)}
       {n \sh \pi n t \mu \, \ch 2 \pi n t \eta}
    \sin 2\pi n i \lam t
    \right)
\\
    &=
    \frac{(q^4 p^{ \lam}      ; p^\mu, q^4)_\infty
          (    p^{ \lam + \mu}; p^\mu, q^4)_\infty
          (q^2 p^{-\lam}      ; p^\mu, q^4)_\infty
          (q^2 p^{-\lam + \mu}; p^\mu, q^4)_\infty}
         {(q^4 p^{-\lam}      ; p^\mu, q^4)_\infty
          (    p^{-\lam + \mu}; p^\mu, q^4)_\infty
          (q^2 p^{ \lam}      ; p^\mu, q^4)_\infty
          (q^2 p^{ \lam + \mu}; p^\mu, q^4)_\infty}
\\
    &=
    \frac{\Gamma_{q^4}\left(\frac{1}{2}+\frac{\lam}{4\eta}\right)
          \Gamma_{q^4}\left(1          -\frac{\lam}{4\eta}\right)}
         {\Gamma_{q^4}\left(\frac{1}{2}-\frac{\lam}{4\eta}\right)
          \Gamma_{q^4}\left(1          +\frac{\lam}{4\eta}\right)}
    \times\\
    &\times
    \prod_{k=1}^\infty
    \frac{\Gamma_{q^4}\left(\frac{1}{2}+\frac{\lam+k\mu}{4\eta}\right)^2
          \Gamma_{q^4}\left(1          +\frac{-\lam+k\mu}{4\eta}\right)
          \Gamma_{q^4}\left(            \frac{-\lam+k\mu}{4\eta}\right)}
         {\Gamma_{q^4}\left(\frac{1}{2}+\frac{-\lam+k\mu}{4\eta}\right)^2
          \Gamma_{q^4}\left(1          +\frac{\lam+k\mu}{4\eta}\right)
          \Gamma_{q^4}\left(            \frac{\lam+k\mu}{4\eta}\right)},
\end{split}
\end{equation}
where $p=e^{-2\pi t}$, $q= e^{-2\pi\eta t}= p^\eta$.
(See \appref{table:func} for definitions of notations. The last equality
is due to \eqref{inf-prod-q-gamma}.) This $\bbbS$ factor was found by Freund
and Zabrodin \cite{fre-zab}.

Comparing \eqref{spec(S)} with \eqref{spec(R)}, we come to the
following conclusion (cf.~\cite{fijkmy}):
\begin{equation}
   S(x) \propto R(\lam;it(1-4\ell\eta)).
\end{equation}

%
%
%
%
\section{Comments and discussions}

\subsection*{}
In this paper we have studied the eigenvectors of transfer matrix of
higher spin generalizations of the eight vertex model by means of
Bethe Ansatz. Apparently Bethe Ansatz for these models seems to be
less powerful compared to Bethe Ans\"atze for the $XXX$ model, the
$XXZ$ model and their higher spin generalizations, since the number of
quasi-particles ($B$ operators) are restricted to $N \ell$. But we
have a discrete parameter $\nu$ instead. In Chapter II, varying this
parameter, we restored all two-particle states which would degenerate
to a singlet and a triplet in the limit, $\eta \to 0$, $t \to \infty$.
Therefore we can expect that Bethe Ansatz for our case gives
as many eigenvectors as that for the rational and trigonometric cases.

\subsection*{}
Developments of the theory of quantum affine algebras in the last
decade provided algebraic tools such as vertex operators and
crystal basis for investigation of the models associated to
trigonometric $R$ matrices. This kind of algebraic method is still
hard to apply to the models associated to elliptic $R$ matrices
because of the lack of knowledge on ``elliptic affine algebras'' which
should be an affinization of the Sklyanin algebra in an appropriate
sense. Foda, Iohara, Jimbo, Kedem, Miwa and Yan \cite{fijkmy} proposed
a candidate of this algebra. In their formulation a relation of the
type $RLL = LLR^\ast$ plays an important role, where $R^\ast$ is
essentially the $S$ marix of two particle states of  the $XYZ$
model. Their argument was based on Smirnov's conjecture, which is
supported by the result of the present paper.

The algebra which Foda et al.~propose is considered to be symmetry of the
$XYZ$ spin chain in a thermodynamic limit. We can also expect that their
algebra is also a symmetry algebra for higher spin models which we
considered in this paper. On the other hand, it is still unknown whether
finite size models could have symmetry of the Sklyanin algebra, since
a reasonable coproduct for the Sklyanin algebra is not yet found.

\subsection*{}
We considered only such excitations with finite number of holes and
finite number of strings which have length $A$, $A\neq 2\ell$, even
in the thermodynamic limit. This is because we wanted to determine the
two particle $S$ matrix. But when we want to calculate thermodynamical
quantities like entropy or specific heat of the model, string
configurations with non-zero hole density and non-zero densities of
$A$-strings ($A \neq 2\ell$) are essential. (See \cite{yang-yang},
\cite{tak-suz}, \cite{kir-resh}.) Extensive thermodynamics
of the $XXX$, $XXZ$ models and their generalization to higher spin
cases has quite interesting features \cite{resh:S}, and are also
related to dilogarithm identities \cite{bab}, \cite{kir-resh}. Further
study of the thermodynamic Bethe Ansatz for higher spin
generalizations of the eight vertex model could give deformation of
the above features in rational and trigonometric models.

\subsection*{Acknowledgements}
The author expresses his gratitude to E.~K.~Sklyanin and J.~Suzuki for
valuable advices and discussions. He also thanks to K.~Hasegawa,
K.~Iohara, M.~Jimbo, A.~N.~Kirillov, V.~Korepin, A.~Kuniba, T.~Miwa,
M.~Noumi, N.~Yu.~Reshetikhin, N.~Slavnov, F.~Smirnov, Y.~Yamada,
A.~Zabrodin for their comments, suggestions and interests.

This work is partly supported by the Grant-in-Aid for Scientific
Research on Priority Areas 231, the Ministry of Education, Science and
Culture.

%
%
%
%
\appendix
\renewcommand{\theequation}{\thesection.\arabic{equation}}
\renewcommand{\thethm}{\thesection.\arabic{thm}}
\section{Review of the Sklyanin algebra}
\setcounter{equation}{0}
\label{skl-alg}
In this appendix we recall several facts on the Sklyanin algebra and
its representations from \cite{skl:alg} and \cite{skl:rep}.
We use notations in \cite{tata} for theta functions:
\begin{equation}
    \theta_{ab}(z;\tau) = \sum_{n\in\Integer}
       \exp\left(
             \pi i \left( \frac{a}{2} + n \right)^2 \tau
           +2\pi i \left( \frac{a}{2} + n \right)
                   \left( \frac{b}{2} + z \right)
           \right),
\label{theta:def}
\end{equation}
where $\tau$ is a complex number such that $\Im(\tau) > 0$.
We denote $t = i/\tau$. The Pauli matrices are defined as usual:
\begin{equation}
    \sig^0 = \begin{pmatrix} 1 & 0  \\  0  & 1 \end{pmatrix},\quad
    \sig^1 = \begin{pmatrix} 0 & 1  \\  1  & 0 \end{pmatrix},\quad
    \sig^2 = \begin{pmatrix} 0 & -i \\  i  & 0 \end{pmatrix},\quad
    \sig^3 = \begin{pmatrix} 1 & 0  \\  0  &-1 \end{pmatrix}.
\end{equation}

The {\em Sklyanin algebra}, $U_{\tau,\eta}(sl(2))$ is generated
by four generators $S^0$, $S^1$, $S^2$, $S^3$, satisfying
the following relations:
\begin{equation}
    R_{12}(\lam-\mu) L_{01}(\lam) L_{02}(\mu) =
    L_{02}(\mu) L_{01}(\lam) R_{12}(\lam-\mu).
\label{RLL}
\end{equation}
Here $\lam$, $\mu$ are complex parameters, the {\em $L$ operator},
$L(\lam)$, is defined by
\begin{equation}
\begin{gathered}
    L(\lam) = \sum_{a=0}^3 W_a^L(\lam) S^a \tensor \sig^a,
\\
    W_0^L(\lam)
    = \frac{\theta_{11}(\lam;\tau)}{\theta_{11}(\eta;\tau)},\qquad
    W_1^L(\lam)
    = \frac{\theta_{10}(\lam;\tau)}{\theta_{10}(\eta;\tau)},
\\
    W_2^L(\lam)
    = \frac{\theta_{00}(\lam;\tau)}{\theta_{00}(\eta;\tau)},\qquad
    W_3^L(\lam)
    = \frac{\theta_{01}(\lam;\tau)}{\theta_{01}(\eta;\tau)},
\end{gathered}
\label{def:L}
\end{equation}
$R(\lam)=R(\lam;it)$ is {\em Baxter's $R$ matrix} defined by
\begin{equation}
    R(\lam) = \sum_{a=0}^3 W_a^R(\lam) \sig^a \tensor \sig^a,\qquad
    W_a^R(\lam) := W_a^L(\lam + \eta).
\label{def:R}
\end{equation}
and indices $\{0,1,2\}$ denote the spaces on which operators act
non-trivially: for example,
\begin{equation*}
    R_{12}(\lam) =
    \sum_{a=0}^3
    W_a^R(\lam) 1 \tensor \sig^a \tensor \sig^a,\qquad
    L_{01}(\lam) =
    \sum_{a=0}^3
    W_a^L(\lam) S^a \tensor \sig^a \tensor 1.
\end{equation*}

The above relation \eqref{RLL} contains $\lam$ and $\mu$ as parameters,
but the commutation relations among $S^a$ ($a=0, \dots, 3$) do not
depend on them:
\begin{equation}
    [S^\alpha,S^0    ]_- =
    -i J_{\alpha,\beta} [S^\beta,S^\gamma]_+, \qquad
    [S^\alpha,S^\beta]_- =
                      i [S^0,    S^\gamma]_+,
\label{comm_rel}
\end{equation}
where $(\alpha, \beta, \gamma)$ stands for any cyclic permutation of
(1,2,3), $[A,B]_\pm = AB\pm BA$, and
$J_{\alpha,\beta}=(W_\alpha^2-W_\beta^2)/(W_\gamma^2-W_0^2)$, i.e.,
\begin{align*}
    J_{12}&= \frac{\theta_{01}(\eta;\tau)^2 \theta_{11}(\eta;\tau)^2}
                 {\theta_{00}(\eta;\tau)^2 \theta_{10}(\eta;\tau)^2},
\\
    J_{23}&= \frac{\theta_{10}(\eta;\tau)^2 \theta_{11}(\eta;\tau)^2}
                 {\theta_{00}(\eta;\tau)^2 \theta_{01}(\eta;\tau)^2},
\\
    J_{31}&=-\frac{\theta_{00}(\eta;\tau)^2 \theta_{11}(\eta;\tau)^2}
                 {\theta_{01}(\eta;\tau)^2 \theta_{10}(\eta;\tau)^2}.
\end{align*}

The {\em spin $\ell$ representation} of the Sklyanin algebra,
$\rho^{\ell}: U_{\tau,\eta}(sl(2))
          \to \End_{\Complex}(\Theta^{4\ell}_{00})$
is defined as follows: The representation space is
\begin{equation}
    \Theta^{4\ell+}_{00} =
    \{f(z) \, |\,
     f(z+1) = f(-z) = f(z), f(z+\tau)=\exp^{-4\ell\pi i(2z+\tau)}f(z) \}.
\end{equation}
It is easy to see that dim$\Theta^{4\ell+}_{00} = 2\ell+1$.
The generators of the algebra act on this space as difference operators:
\begin{equation}
    (\rho^\ell(S^a) f)(z) =
    \frac{s_a(z-\ell\eta)f(z+\eta)-s_a(-z-\ell\eta)f(z-\eta)}
         {\theta_{11}(2z;\tau)},
\end{equation}
where
\begin{alignat*}{2}
    s_0(z) &=  \theta_{11}(\eta;\tau) \theta_{11}(2z;\tau),\qquad&
    s_1(z) &=  \theta_{10}(\eta;\tau) \theta_{10}(2z;\tau),\\
    s_2(z) &= i\theta_{00}(\eta;\tau) \theta_{00}(2z;\tau),\qquad&
    s_3(z) &=  \theta_{01}(\eta;\tau) \theta_{01}(2z;\tau).
\end{alignat*}
These representations reduce to the usual spin $\ell$ representations of
$U(sl(2))$ for $J_{\alpha\beta} \to 0$ ($\eta \to 0$).
In particular, in the case $\ell= 1/2$, $S^a$  are expressed by
the Pauli matrices $\sigma^a$ as follows:
Take
$(\theta_{00}(2z;2\tau)-\theta_{10}(2z;2\tau),
  \theta_{00}(2z;2\tau)+\theta_{10}(2z;2\tau))$ as a basis of
$\Theta^{2+}_{00}$. With respect to this basis $S^a$ have matrix forms
\begin{equation}
    \rho^{1/2}(S^a) =
          2\frac{\theta_{00}(\eta;\tau) \theta_{01}(\eta;\tau)
                 \theta_{10}(\eta;\tau) \theta_{11}(\eta;\tau)}
                {\theta_{00}(0;\tau) \theta_{01}(0;\tau)
                 \theta_{10}(0;\tau) \hfill}
          \sigma^a.
\label{rep:pauli}
\end{equation}
Since the relations \eqref{comm_rel} are homogeneous, an overall
constant factor in a representation is not essential.

There are involutive automorphisms of the Sklyanin algebra
$U_{\tau,\eta}(sl(2))$ found by Sklyanin \cite{skl:rep}:
for $a=1,2,3$,
\begin{equation}
    X_a: (S^0, S^a,  S^b,  S^c) \mapsto
         (S^0, S^a, -S^b, -S^c),
\label{def:X-a}
\end{equation}
where $(a,b,c)$ is a cyclic permutation of $(1,2,3)$.
Combining these operators with $\rho^{\ell}$, we obtain another
representation $\rho^\ell \circ X_a$ of $U_{\tau,\eta}(sl(2))$,
but there is a unitary operator $U_a$ intertwining $\rho^{\ell}$
and $\rho^\ell \circ X_a$ \cite{skl:rep}:
\begin{equation*}
\begin{alignedat}2
    U_1: &\Theta^{4\ell+}_{00} \owns f(z) \mapsto &
         &(U_1 f)(z) = e^{\pi i \ell} f\left(z + \frac{1}{2} \right),
\\
    U_3: &\Theta^{4\ell+}_{00} \owns f(z) \mapsto &
         &(U_3 f)(z) = e^{\pi i \ell} e^{\pi i \ell (4z+\tau)}
                       f\left(z + \frac{\tau}{2} \right),
\end{alignedat}
\end{equation*}
and $U_2= U_3 U_1$. Direct calculations show that
$X_a(\rho^\ell(S^b)) = U_a^{-1} \rho^\ell(S^b) U_a$.
Operators $U_a$ satisfy the relations: $U_a^2 = (-1)^{2\ell}$,
$U_a U_b = (-1)^{2\ell} U_b U_a = U_c$.

Baxter's $R$ matrix \eqref{def:R} is a $4\times4$ matrix proportional
to that in \cite{takh-fad:xyz}:
\begin{equation}
    R(\lam;it)=
    \begin{pmatrix}
         a(\lam) & 0       & 0       & d(\lam) \\
         0       & c(\lam) & b(\lam) & 0       \\
         0       & b(\lam) & c(\lam) & 0       \\
         d(\lam) & 0       & 0       & a(\lam) \\
    \end{pmatrix},
\end{equation}
where functions $a$, $b$, $c$, $d$ are defined by
\begin{align*}
    a(\lam)
    &= C_1\theta_{01}(2it\eta;2it)
          \theta_{01}( it\lam;2it)
          \theta_{11}( it\lam + 2it\eta;2it),
\\
    b(\lam)
    &= C_1\theta_{11}(2it\eta;2it)
          \theta_{01}( it\lam;2it)
          \theta_{01}( it\lam + 2it\eta;2it),
\\
    c(\lam)
    &= C_1\theta_{01}(2it\eta;2it)
          \theta_{11}( it\lam;2it)
          \theta_{01}( it\lam + 2it\eta;2it),
\\
    d(\lam)
    &= C_1\theta_{11}(2it\eta;2it)
          \theta_{11}( it\lam;2it)
          \theta_{11}( it\lam + 2it\eta;2it),
\\
    C_1 &= \frac{- 2 \exp(- \pi t \lam (\lam + 2\eta))}
            {\theta_{01}(0;2it)
             \theta_{01}(2it\eta;2it)
             \theta_{11}(2it\eta;2it)}.
\end{align*}
Obviously eigenvalues of this matrix are
\begin{equation}
\begin{alignedat}2
    a(\lam) + d(\lam) &= C_2
    \frac{\theta_{00}\left(\frac{it\lam}{2} - it\eta;it\right)}
         {\theta_{00}\left(\frac{it\lam}{2} + it\eta;it\right)},
\quad&
    a(\lam) - d(\lam) &= C_2
    \frac{\theta_{01}\left(\frac{it\lam}{2} - it\eta;it\right)}
         {\theta_{01}\left(\frac{it\lam}{2} + it\eta;it\right)},
\\
    b(\lam) + c(\lam) &= C_2
    \frac{\theta_{10}\left(\frac{it\lam}{2} - it\eta;it\right)}
         {\theta_{10}\left(\frac{it\lam}{2} + it\eta;it\right)},
\quad&
    b(\lam) - c(\lam) &= C_2
    \frac{\theta_{11}\left(\frac{it\lam}{2} - it\eta;it\right)}
         {\theta_{11}\left(\frac{it\lam}{2} + it\eta;it\right)},
\end{alignedat}
\label{spec(R)}
\end{equation}
where
\begin{equation*}
\begin{split}
    C_2 =& 2 e^{-\pi t \lam (\lam + 2\eta)} \times\\
    \times&\frac{\theta_{00}\left(\frac{it\lam}{2} + it\eta;it\right)
          \theta_{01}\left(\frac{it\lam}{2} + it\eta;it\right)
          \theta_{10}\left(\frac{it\lam}{2} + it\eta;it\right)
          \theta_{11}\left(\frac{it\lam}{2} + it\eta;it\right)}
         {\theta_{00}(it\eta;it)
          \theta_{01}(it\eta;it)
          \theta_{10}(it\eta;it)
          \theta_{11}(it\eta;it)}.
\end{split}
\end{equation*}

%
%
%
%
\section{Proof of the sum rule}
\setcounter{equation}{0}
\label{sumrule:pf}
We prove here the sum rule of $\lam_j$'s, \thmref{sumrule:thm}.
See \secref{aBA} for notations.

Let us introduce a determinant $t^r(\lam)$ of a $(r-1)\times (r-1)$
matrix, elements of which are defined by:
\begin{enumerate}
\renewcommand{\labelenumi}{(\arabic{enumi})}
\item $(j,j+1)$-elements $= h(\lam+2(j-\ell)\eta)$;
\item $(j,j)$-elements $= t(\lam+2j\eta)$;
\item $(j,j-1)$-elements $=h(\lam+2(j+\ell)\eta)$;
\item other elements are 0,
\end{enumerate}
where $t(\lam)$ is the eigenvalue of the transfer matrix $T(\lam)$
on the Bethe vector $\Psi_\nu(\lam_1, \dots, \lam_M)$ and
$h(z) = (2\theta_{11}(z))^N$.
(This determinant is related to a fused model.)
Since $t(\lam)$ is an entire function of $\lam$ (recall that the
transfer matrix $T(\lam)$ itself is an entire function of $\lam$),
$t^r(\lam)$ is an entire function of $\lam$. (In other words, the
analyticity is a consequence of the Bethe equations as noted in
\secref{aBA}.) Our third assumption (see \thmref{sumrule:thm}) is

\medskip
iii) $t^r(\lam)$ is not identically zero.

\medskip

Let us define $\tilde t^r(\lam)$ by
\begin{equation}
    \tilde t^r(\lam) :=
    Q(\lam + 2\eta) Q(\lam + 4\eta) \dots Q(\lam + 2(r-1)\eta)
    t^r(\lam).
\label{def:tt-r}
\end{equation}
Then because of \eqref{Bethe-ev:Q}
\begin{equation*}
    t(\lam):= h(\lam+2\ell\eta) \frac{Q(\lam-2\eta)}{Q(\lam)}
            + h(\lam-2\ell\eta) \frac{Q(\lam+2\eta)}{Q(\lam)},
\end{equation*}
function $\tilde t^r(\lam)$ is expressed as a determinant of a matrix
such that
\begin{enumerate}
\renewcommand{\labelenumi}{(\arabic{enumi})}
\item $(j,j+1)$-element $= a_+(\lam+2j\eta)$;
\item $(j,j)$-element $= a_-(\lam+2j\eta) + a_+(\lam+2j\eta)$;
\item $(j,j-1)$-element $= a_-(\lam+2j\eta)$,
\item other elements are 0.
\end{enumerate}
where
\begin{align*}
    a_-(\lam) &= h(\lam + 2\ell\eta) Q(\lam - 2\eta),\\
    a_+(\lam) &= h(\lam - 2\ell\eta) Q(\lam + 2\eta).
\end{align*}
This determinant can be easily expanded, the result being
\begin{equation}
\begin{split}
    \tilde t^r(\lam) =&
    \sum_{j=1}^r
    a_-(\lam+2\eta) \dots a_-(\lam+2(j-1)\eta)
    a_+(\lam+2j\eta) \dots a_+(\lam+2(r-1)\eta)\\
    =&
    h(\lam+2(\ell+1)\eta) \dots h(\lam+2(r-\ell-1)\eta) Q(\lam)
    (f_0(\lam) + \dots + f_{r-1}(\lam)),
\end{split}
\end{equation}
where
\begin{equation*}
    f_k(\lam) = \prod_{j=1}^{2\ell} h(\lam+2(k-\ell+j)\eta)
                \prod_{j=1}^{k-1} Q(\lam+2j\eta)
                \prod_{j=k+2}^r Q(\lam+2j\eta).
\end{equation*}
By the definition \eqref{def:Q}, $Q(\lam)$ has automorphic property:
\begin{equation}
\begin{aligned}
    Q(\lam+1)    &= (-1)^{N\ell - \nu} Q(\lam),\\
    Q(\lam+\tau) &= e^{-\pi i N\ell(1+\tau+2\lam) - \pi i \tau \nu
                    + 2 \pi i \sum_{j=1}^M \lam_j} Q(\lam).
\end{aligned}
\label{Q:auto}
\end{equation}
Hence $f_k(\lam+2\eta) = f_{k+1}(\lam)$ ($f_r(\lam)=f_0(\lam)$) and
$F(\lam) =f_0(\lam) + \dots + f_{r-1}(\lam)$ has a period $2\eta$:
$F(\lam+2\eta) = F(\lam)$.

Now we proceed in four steps.

\subsection*{Step1}
First we show that
$Q(\lam)F(\lam)/Q(\lam+2\eta) \dots Q(\lam+2(r-1)\eta)$ is an entire
function of $\lam$. Since
\begin{equation}
\begin{split}
    t^r(\lam) &= \frac{\tilde t^r(\lam)}
                      {Q(\lam+2\eta) \dots Q(\lam+2(r-1)\eta)}\\
    &=
    h(\lam+2(\ell+1)\eta) \dots h(\lam+2(r-\ell-1)\eta)
    \frac{Q(\lam)F(\lam)}
         {Q(\lam+2\eta) \dots Q(\lam+2(r-1)\eta)}
\end{split}
\end{equation}
is an entire function of $\lam$, we have only to show that any zero of
the denominator is not a zero of
$h(\lam+2(\ell+1)\eta) \dots h(\lam+2(r-\ell-1)\eta)$.
Zeros of $h(\lam)$ is 0 $\mod \Integer+\Integer\tau$. Hence the last
statement is true if assumption i) of \thmref{sumrule:thm} is
fulfilled.

\subsection*{Step2}
We show that $F(\lam)/Q(\lam+2\eta) \dots Q(\lam+2(r-1)\eta)$ is
an entire function of $\lam$.

As a consequence of Step1, we know that only possible poles of
$F(\lam)/Q(\lam+2\eta) \dots Q(\lam+2(r-1)\eta)$ exist at zeros of
$Q(\lam)$. Suppose $\lam_j$ is a pole of
$F(\lam)/Q(\lam+2\eta) \dots Q(\lam+2(r-1)\eta)$.
Then
\begin{equation}
    \ord_{\lam_j} F(\lam) <
    \ord_{\lam_j} (Q(\lam+2\eta) \dots Q(\lam+2(r-1)\eta)) \leqq
    \ord_{\lam_j} F(\lam) + \ord_{\lam_j} Q(\lam).
\label{step2.1}
\end{equation}
Here $\ord_{\lam_j}$ is the order of zero at $\lam_j$.
Assumption ii) of \thmref{sumrule:thm} says that there is an integer
$a$ such that $\ord_{\lam_j+2a\eta} Q(\lam) = 0$. On the other hand
periodicity $F(\lam+2\eta) = F(\lam)$ implies that
$\ord_{\lam_j+2a\eta} F(\lam)=\ord_{\lam_j} F(\lam)$. Therefore
\begin{equation}
\begin{split}
     & \ord_{\lam_j+2a\eta} (Q(\lam+2\eta) \dots Q(\lam+2(r-1)\eta))\\
    =& \ord_{\lam_j} (Q(\lam+2(a+1)\eta) \dots Q(\lam+2(r+a-1)\eta))\\
    =& \ord_{\lam_j} (Q(\lam+2 a   \eta) \dots Q(\lam+2(r+a-1)\eta))\\
    =& \ord_{\lam_j} (Q(\lam+2     \eta) \dots Q(\lam+2(r  -1)\eta))\\
    >& \ord_{\lam_j} F(\lam) = \ord_{\lam_j+2a\eta} F(\lam)\\
    =& \ord_{\lam_j+2a\eta} Q(\lam) F(\lam)
\end{split}
\label{step2.2}
\end{equation}
Here we used \eqref{Q:auto} and \eqref{step2.1}.
This inequality means that
$Q(\lam) F(\lam)/(Q(\lam+2\eta) \dots Q(\lam+2(r-1)\eta))$ has a pole
at $\lam_j + 2a\eta$. This contradicts Step 1.

\subsection*{Step3}
Now we show that even
$F(\lam)/(Q(\lam) Q(\lam+2\eta) \dots Q(\lam+2(r-1)\eta))$ is
an entire function of $\lam$. It follows from Step 2 that
for any $j=0, 1, \dots, r-1$
\begin{equation*}
    \frac{Q(\lam+2j\eta) F(\lam)}
         {Q(\lam) Q(\lam+2\eta) \dots Q(\lam + 2(r-1)\eta)}
\end{equation*}
is an entire function. Suppose
$F(\lam)/(Q(\lam) Q(\lam+2\eta) \dots Q(\lam+2(r-1)\eta))$ has a pole
at $\lam_0$. Then $\lam_0$ should be a zero of $Q(\lam+2j\eta)$,
$j=0, \dots, r-1$. Taking \eqref{Q:auto} into account, this
contradicts assumption ii).

\subsection*{Step4}
We have shown that $F(\lam)/G(\lam)$ is an entire function where
$G(\lam) = Q(\lam) Q(\lam+2\eta) \dots Q(\lam+2(r-1)\eta)$.
Using \eqref{Q:auto} and
\begin{align*}
    h(\lam+1)    &= (-1)^N h(\lam),\\
    h(\lam+\tau) &= e^{- \pi i N(1+\tau)- 2\pi i \lam} h(\lam),
\end{align*}
we obtain
\begin{align}
    \frac{F(\lam+1)}{G(\lam+1)} &= \frac{F(\lam)}{G(\lam)},
\label{F/G:period}
\\
    \frac{F(\lam+\tau)}{G(\lam+\tau)} &=
    e^{2\pi i (\nu\tau - 2 \sum_{j=1}^M \lam_j)}
    \frac{F(\lam)}{G(\lam)}.
\label{F/G:tau}
\end{align}
Holomorphy of $F/G$ and periodicity \eqref{F/G:period} make it
possible to expand $F/G$ into a Fourier series:
\begin{equation*}
    (F/G)(\lam) = \sum_{n\in\Integer} (F/G)_n e^{2\pi i n \lam}.
\end{equation*}
Substituting $\lam+\tau$ into this expansion and comparing with
\eqref{F/G:tau}, we find that each coefficient should satisfy
\begin{equation*}
    (F/G)_n =
    (F/G)_n e^{2\pi i ((\nu - n)\tau - 2\sum_{j=1}^M \lam_j)}.
\end{equation*}
Since $\Im\tau >0$, there exists only one $n=n_2$ such that
$(F/G)_{n_2}\neq0$ and it satisfies
$(\nu - n_2)\tau - 2\sum_{j=1}^M \lam_j =: - n_0 \in \Integer$.
Putting $n_1 = \nu - n_2$, we have
$2\sum_{j=1}^M \lam_j = n_0 + n_1\tau$.

It follows from the above argument that $t^r(\lam)$ has the following
form with a suitable integer $n$:
\begin{equation}
    t^r(\lam) = \text{ const. }e^{2\pi i n \lam}
    h(\lam+2(\ell+1)\eta) \dots h(\lam+2(r-\ell-1)\eta) Q(\lam)^2.
\end{equation}

%
%
%
%
\section{Table of useful functions}
\setcounter{equation}{0}
\label{table:func}
Here we collect properties of functions used in Chapter II.

\subsection*{Logarithm of quotient of theta functions}

A function $\Phi$ defined by \eqref{def:Phi},
\begin{equation*}
    \Phi(x;i\mu t) = \frac{1}{i}
                   \log \frac{\theta_{11}(x + i\mu t ;it)}
                             {\theta_{11}(x - i\mu t ;it)} + \pi,
\end{equation*}
has the following Fourier expansion if $0 < \mu < 1/2$:
\begin{equation}
    \Phi(x;i\mu t)=
    -2\pi x
    - 2 \sum_{n=1}^\infty
    \frac{\sh \pi n (1-2\mu) t}{n \sh \pi n t}
    \sin 2\pi n x.
\label{Fourier:Phi}
\end{equation}
Hence
\begin{equation}
    \frac{d}{dx}\Phi(x;i\mu t) =
    -2\pi\left(1+ 2 \sum_{n=1}^\infty
                \frac{\sh \pi n (1-2\mu) t}{\sh \pi n t}
                \cos 2\pi n x \right).
\label{Fourier:Phi'}
\end{equation}

A function $\Psi$ defined by \eqref{def:Psi},
\begin{equation*}
    \Psi(x;i\mu t) = \frac{1}{i}
                   \log \frac{\theta_{01}(x + i\mu t ;it)}
                             {\theta_{01}(x - i\mu t ;it)},
\end{equation*}
has the following Fourier expansion if $0 < \mu < 1/2$:
\begin{equation}
    \Psi(x;i\mu t)=
    2 \sum_{n=1}^\infty
    \frac{\sh 2 \pi n \mu t}{n \sh \pi n t}
    \sin 2\pi n x.
\label{Fourier:Psi}
\end{equation}
Hence
\begin{equation}
    \frac{d}{dx}\Psi(x;i\mu t) =
           4\pi \sum_{n=1}^\infty
                \frac{\sh 2 \pi n \mu t}{\sh \pi n t}
                \cos 2\pi n x .
\label{Fourier:Psi'}
\end{equation}

\begin{lem}
\label{positivity}
For $0<a<b$, a series
\begin{equation*}
    \sum_{n\in\Integer} \frac{\sh\pi n a}{\sh\pi n b} e^{2\pi i nx}
\end{equation*}
is positive for $x\in\Real$. Here the term $n=0$ is understood as
$a/b$.
\end{lem}

\begin{pf}
Define a function $f(y;x)$ by
\begin{equation*}
    f(y;x) = \frac{\sh\pi y a}{\sh\pi y b} e^{2\pi i yx},
\end{equation*}
$f(0;x) = a/b$. The Fourier transformation of this function is
\begin{equation}
\begin{split}
    \hat f(\xi;x)
    &= \int_{-\infty}^\infty f(y;x) e^{-2\pi i y \xi}\, dy\\
    &= \frac{2\sin(a\pi/b)}{b}
       \frac{e^{-2\pi|x-\xi|/b}}
       {\left(e^{-2\pi|x-\xi|/b} + \cos (a\pi/b)\right)^2 +
        \sin^2(a\pi/b)} >0
\end{split}
\end{equation}
By Poisson's summation formula we have
\begin{equation*}
    \sum_{n\in\Integer} \frac{\sh\pi n a}{\sh\pi n b} e^{2\pi i nx}
    =
    \sum_{n\in\Integer} f(n;x)
    =
    \sum_{n\in\Integer} \hat f(n;x) >0.
\end{equation*}
This proves the lemma.
\end{pf}

For example $d/dx (\Phi(x;i\mu t)) < 0$, because of \lemref{positivity}.

\subsection*{$q$-$\Gamma$ function}

A $q$-analogue of the $\Gamma$ function is defined by
\begin{equation}
    \Gamma_q(x) = \frac{(q;q)_\infty}{(q^x;q)_\infty} (1-q)^{1-x},
\label{def:q-gamma}
\end{equation}
where $(x;q)_\infty = \prod_{n=0}^\infty (1-xq^n)$.

Double infinite product $(x;q_1, q_2)_\infty$ is defined by
$ (x;q_1, q_2)_\infty
= \prod_{n_1,n_2 = 0}^\infty (1-xq_1^{n_1}q_2^{n_2})$.
If $x+y=z+w$, the following relation holds:
\begin{equation}
    \frac{(q^x; q^a, q)_\infty (q^y; q^a, q)_\infty}
         {(q^z; q^a, q)_\infty (q^w; q^a, q)_\infty}
    =
    \prod_{n=1}^\infty
    \frac{\Gamma_q(z+an) \Gamma_q(w+an)}
         {\Gamma_q(x+an) \Gamma_q(y+an)}
\label{inf-prod-q-gamma}
\end{equation}

%
%
%
%

%
\end{document}